\documentclass{tac}[11pt]

\usepackage[all,cmtip,2cell]{xy}
\usepackage{mathtools}
\usepackage{mathrsfs}
\usepackage{amsmath}
\usepackage{hyperref}
\usepackage{url}
\usepackage{mathrsfs}
\usepackage{amsfonts}
\usepackage{amssymb}
\UseTwocells
\numberwithin{equation}{section}

\newcommand{\E}{\mathscr E}
\newcommand{\A}{\mathscr A}
\newcommand{\B}{\mathscr B}
\newcommand{\D}{\mathbb D}

\newcommand{\cat}{\mathbf{Cat}}

\newcommand{\set}{\mathbf{Set}}
\newcommand{\dblset}{\mathbb{S}\mathbf{et}}
\newcommand{\Span}{\mathbb S\mathbf{pan}}

\newcommand{\Rel}{\mathbb R\mathbf{el}}

\newcommand{\elt}{\mathbf{Elt}}

\newcommand{\dbl}{\mathbf{Dbl}}

\newcommand{\prof}{\mathbb P\mathbf{rof}}
\newcommand{\alg}{\mathbf{Alg}}

\newcommand{\copt}{\mathbf{Copt}}
\newcommand{\Sub}{\mathrm{Sub}}
\newcommand{\dopf}{\mathbf{DOpf}}

\newcommand{\src}{\mathrm{src}}
\newcommand{\tgt}{\mathrm{tgt}}
\newcommand{\hto}{\slashedrightarrow}

\makeatletter
\def\slashedarrowfill@#1#2#3#4#5{%
  $\m@th\thickmuskip0mu\medmuskip\thickmuskip\thinmuskip\thickmuskip
  \relax#5#1\mkern-7mu%
  \cleaders\hbox{$#5\mkern-2mu#2\mkern-2mu$}\hfill
  \mathclap{#3}\mathclap{#2}%
  \cleaders\hbox{$#5\mkern-2mu#2\mkern-2mu$}\hfill
  \mkern-7mu#4$%
}
\def\rightslashedarrowfill@{%
  \slashedarrowfill@\relbar\relbar\mapstochar\rightarrow}
\newcommand\xslashedrightarrow[2][]{%
  \ext@arrow 0055{\rightslashedarrowfill@}{#1}{#2}}
\makeatother
\def\slashedrightarrow{\xslashedrightarrow{}}

\title{Double Categories of Relations}

\author{Michael Lambert}

\copyrightyear{2022}

\address{Mathematics Department, University of Massachusetts-Boston \\100 William T. Morrissey Boulevard\\
Boston, MA 02125 \\
}

\eaddress{michael.lambert@umb.edu, michael.james.lambert@gmail.com
}

\begin{document}

\maketitle

\begin{abstract}
    A `double category of relations' is defined in this paper as a cartesian equipment in which every object is suitably discrete. The main result is a characterization theorem that a `double category of relations' is equivalent to a double category of relations on a regular category when it has strong and monic tabulators and a double-categorical subobject comprehension scheme. This result is based in part on the recent characterization of double categories of spans due to Aleiferi. The overall development can be viewed as a double-categorical version of that of the notion of a ``functionally complete bicategory of relations" or a ``tabular allegory."
\end{abstract}

\tableofcontents

\section{Introduction}

`Bicategories of relations' were introduced in \cite{CW}. This paper aims to give a double-categorical version of this development, leading to a definition of a `double category of relations' and a characterization theorem specifying conditions under which a given double category is a double category of relations on a regular category. More specifically, a `\textbf{bicategory of relations}' is a cartesian bicategory in which every object is suitably discrete. The characterization theorem of \cite{CW} says that every functionally complete `bicategory of relations' is equivalent to one of the form $\mathfrak{Rel}(\E)$. ``Functionally complete" means that every arrow has a so-called ``tabulation." It will be seen in this paper that the horizontal bicategory of any cartesian equipment is a cartesian bicategory, meaning that it makes sense to define a `double category of relations' as a cartesian equipment satisfying the appropriate discrete condition. Just as every functionally complete `bicategory of relations' is of the form $\mathfrak{Rel}(\E)$, it will be seen that every suitably functionally complete `double category of relations' is a double category of the form $\Rel(\E)$ for some regular category $\E$.

This paper has a narrative purpose, intended to explain the naturalness of the conditions involved in the definition of a functionally complete `double category of relations'. First of all, `bicategories of relations' are \emph{cartesian} bicategories. Thus, the starting point is the recent definition \cite{aleiferi} of a ``cartesian double category" and the accompanying characterization of cartesian double categories of the form $\Span(\E)$ for some finitely-complete $\E$. A double category $\D$ is equivalent to one of the form $\Span(\E)$ for some finitely-complete $\E$ if, and only if, $\D$ is a unit-pure cartesian equipment with strong tabulators and certain Eilenberg-Moore objects. What distinguishes the present approach from that of \cite{aleiferi} is the emphasis on tabulators, owing to their centrality for `bicategories of relations'. In an arbitrary double category, tabulators are given by a right adjoint $\top\colon \D_1\to \D_0$ to the external identity $y\colon \D_0\to \D_1$ \cite{GP2}. That tabulators exist and are ``strong'' and ``monic'' is the definition of a ``functionally complete'' cartesian equipment given in \S \ref{section:RelsTabulatorsAreStrongandMonic}.

The approach to proving the characterization theorem has two parts. In \S \ref{section:OplaxLaxAdjunctionConditions} and \S \ref{section:ExistenceOplaxLaxAdjEquiv} conditions are developed under which a double category $\D$ is equivalent to $\Rel(\D_0)$. This follows \cite{Niefield} that gives conditions under which a double category $\D$ admits an oplax/lax adjunction to $\Span(\D_0)$. The corresponding adjointness in the relations case reveals conditions under which it is a strong equivalence. These are presented in Theorems \ref{theorem:PreliminaryCharacterizationTheorem} and \ref{theorem:PreliminaryCharacterizationTheoremRedux}. Secondly, \S \ref{section:SiteForRelations} develops conditions under which $\D_0$ is a regular category, meaning that $\Rel(\D_0)$ makes sense in the first place. Familiar properties and constructions enter into this development. For example, a Modular Law is developed in \S \ref{equation:ModularLaw} which is used in the proof of Theorem \ref{theorem:D0isRegular} showing that $\D_0$ admits a suitably rich factorization system for forming $\Rel(\D_0)$. Additionally, ``unit-pure" is required to understand the relationship between monic arrows in $\D_0$ and inclusions in the equipment structure on $\D$. 

However, these are embraced by a couple of innovations. It is shown in Proposition \ref{prop:HorizontalBicatIsCartesianBicat} that the horizontal bicategory of any cartesian double category is a cartesian bicategory. This leads to the definition of a `double category of relations' in \S \ref{section:DoubleCategoriesOfRelations} as a cartesian double category satisfying the analogue of the Frobenius Law from \cite{CW}. As a consequence, in any `double category of relations' the Modularity Law will hold, since this is true in its horizontal bicategory. Secondly, the exactness and functoriality conditions appearing in \S \ref{section:ExistenceOplaxLaxAdjEquiv} are embraced by one concerning the existence of a ``subobject comprehension scheme.'' These are discussed in \S \ref{section:Comprehenseionscheme}. Briefly, interpreting tabulators as generalized elements constructions, the exactness and functoriality conditions are subsumed by asking that tabulators have a pseudo-inverse, or ``fibers'' construction, leading to a certain equivalence of categories.

The main result of the paper, namely, Theorem \ref{theorem:MainCharacterizationTheorem} is that these two conditions are sufficient for $\D$ to be of the form $\Rel(\E)$. This is 

\theorem
    If $\D$ is a `double category of relations' with a subobject comprehension scheme, then $\D_0$ is regular and $1\colon \D_0\to\D_0$ extends to an adjoint equivalence $\Rel(\E)\simeq \D$.
\endtheorem

Once this result is completely proved, some consequences and extensions of the theory are discussed in the remainder of \S \ref{section:CharacterizationTheorem}. It is shown that double categories of relations are regular equipments in the sense of \cite{Schultz}. A concluding section \S \ref{section:DivisionPowers} presents definitions of ``division'' and ``powers,'' giving analogues of those in \cite{FS}. 

The theory presented here has many potential applications. Some of these are discussed in a prospectus in \S \ref{section:Prospectus}. The presentation is only in outline since this paper is a theoretical one. Briefly, several of these applications have to do with using `double categories of relations' as a meeting place, or single axiomatic forum, for ``functional'' and ``relational'' approaches to various topics in (applied) category theory. For example, the ``functional ologs'' of \cite{SpivakKent} have their counterpart in the ``relational ologs'' of \cite{patterson}. A first pass on combining these appears in \cite{lamberttoposblog}. Traditional algebraic theories have their relational version in the ``relational'' and perhaps ``partial'' theories of for example \cite{BFS} and \cite{LLNSob}. Other potential applications concern interpretation of various type theories and logics \cite{BJ}, double categorical models of databases \cite{Codd}, \cite{RW} and of asynchronous communication in distributed systems \cite{SelingerAsynchronous}, and finally with classifications of monoidal bifibrations \cite{FramedBicats}. The last important point worth mentioning is that the ``subobject classification schemes'' discussed here ought to be a central ingredient in some hypothetical notion of a ``double topos'' which is yet to be fully formulated, but should generalize the notion of a 2-topos \cite{Weber}.

Details on double categories can be found in \cite{GP} and \cite{GP2}. Throughout double categories will always be ``pseudo," that is, pseudo-categories in $\cat$ viewed as a 2-category. Conventions and notation for the most part are adopted from those in \cite{FrameworkGenMulticats} and \cite{Schultz}. Many of these are summarized in a previous paper \cite{discdblefibs}. Double categories are denoted in blackboard font such as $\D$. The objects of the category $\mathbb D_1$ are called ``proarrows'' and denoted with a slashed arrow $p\colon A\hto B$. External composition of proarrows will always be written in diagrammatic order. So, $p\otimes q$ for proarrows $p$ and $q$ means ``first $p$ then $q$.''  An \textbf{oplax/lax adjunction} between double categories is a conjoint pair of arrows in the strict double category of lax and oplax functors \cite{GP2}. An oplax/lax adjunction is \textbf{strong} if both functors are pseudo. The work in this paper is meant to generalize the account of `bicategories of relations', but some tools from the theory of allegories will enter into the discussion. Recall that an \textbf{allegory} \cite{FS} is a locally-ordered 2-category equipped with local products and an anti-involution $(-)^\circ$ satisfying a Modularity Law. Each ``unitary" and ``tabular" allegory is equivalent to $\mathbf{Rel}(\E)$ for some $\E$. In this sense, tabular allegories are another axiomatization of a calculus of relations. Allegories and `bicategories of relations' are closely related \cite{PetrusAndFrank}, so it is partly a matter of taste which is chosen as basic. The present use of double-categorical analogues of tools from each account will hopefully exhibit this closeness in practice. 

\subsection{acknowledgements}
    A draft of the results in this paper was presented at ACT 2021. An improved version was later presented at the Dalhousie University AtCat seminar. The results and paper itself evolved considerably over several drafts and would not have achieved their present state without feedback, questions, comments and encouragement of Eva Aleiferi, Bryce Clarke, Geoff Cruttwell, Fosco Loregian, Chad Nester, Bob Par\'e, Dorette Pronk, Martin Szyld and Richard Wood. Special thanks are due to Geoff Cruttwell under whose direct supervision the final version of this paper was completed during the author's postdoc at Mount Allison University. Thanks finally to the editors at \emph{TAC} and to the referee for a very careful reading of the paper.

\section{Double Category Structure of Relations}
\label{section:StructureOfRelations}

A calculus of relations is supported by any regular category $\E$. This leads to the formation of a double category $\Rel(\E)$ whose salient structures are analyzed in this section. In particular, $\Rel(\E)$ is cartesian as a double category and an equipment, meaning that it has all extensions and restrictions, as explained below. These structures provide minimal expectations of any `double category of relations'. They should also satisfy a modular law.

Recall \cite[\S A1.3]{Elephant}, \cite[\S I,1.5]{FS} that the \textbf{image} of a morphism $f\colon A\to B$ is the smallest subobject of $B$ through which $f$ factors, if it exists. Such a factorization $f=me$ is an \textbf{image factorization}. A \textbf{cover} is a morphism $f\colon A\to B$ whose image is all of $B$. Denote covers using `$\twoheadrightarrow$' and monics by `$\rightarrowtail$'. 

\definition \label{define:RegularCat}
    A cartesian category $\E$ is \textbf{regular} if 
        \begin{enumerate}
            \item every morphism has an image factorization;
            \item covers are pullback-stable.
        \end{enumerate}
\enddefinition

These conditions imply the familiar ones, namely, that every kernel has a coequalizer and that regular epimorphisms are pullback-stable \cite[\S I,1.566]{FS}. A \textbf{relation} is a monic arrow $R\rightarrowtail A\times B$. Denote these by `$R\colon A\hto B$'. For any regular category $\E$, take $\Rel(\E)$ to denote the double category of relations in $\E$. Its ordinary underlying category is $\Rel(\E)_0 = \E$. Its proarrows are relations and its cells $\theta$, as at left below, are morphisms $\theta\colon R\to S$ making the square on the right commute:
    $$\xymatrix{
        \ar@{}[dr]|{\theta} A \ar[d]_f \ar[r]^R|-@{|} & B \ar[d]^g & & & R \ar[d]_\theta \ar[r] & A\times B  \ar[d]^{f\times g} \\
        C \ar[r]_S|-@{|} & D & & & S \ar[r] & C\times D
    }$$
Take $y\colon \E\to\Rel(\E)$ to denote the functor sending an object $X$ to the span consisting of identity arrows on $X$. This is fully faithful by construction. The external source and target functors $\src,\tgt\colon \Rel(\E) \rightrightarrows \E$ are given by taking a cell $\theta$ as above to $f$ and to $g$, respectively. What will be external composition is given by pullback in $\E$ and taking images. That is, given two relations $R\to A\times B$ and $S\to B\times C$ take the pullback $R\times_BS$ in $\E$ and define the composite $R\otimes S$ to be given by the image factorization of $R\times_BS\to A\times C$; the arrow assignment is then induced by minimality of image factorizations. Up-to-iso associativity follows by pullback-stability of covers. With these constructions $\Rel(\E)$ is a double category. It is rich with structure. In particular, $\Rel(\E)$ is an ``equipment,'' in the sense immediately below. The importance of this condition is that it seems the minimal one enabling various constructions in formal category theory \cite{roalddoubledim}, \cite{roaldaugmented}.

\definition[\S 4 \cite{FramedBicats}]
    A double category $\D$ is an \textbf{equipment} if the source-target projection functor $\langle\src,\tgt\rangle\colon \D_1\to\D_0\times\D_0$ is a bifibration.
\enddefinition

In the literature, equipments have also been called ``framed bicategories" \cite{FramedBicats} and ``fibrant double categories'' \cite{aleiferi}. The present choice of ``equipment'' goes back to the notion of a ``2-category equipped with proarrows'' \cite{WoodProI}, \cite{WoodProII}. In any case, a bit of definitional exegesis is in order. A cell $\theta$ in $\D$ is \textbf{cartesian} if, and only if, it is a cartesian arrow for the functor $\langle \src,\tgt\rangle\colon \D_1\to \D_0\times\D_0$, that is, if given any other cell $\delta$
    $$\xymatrix{
        \ar@{}[dr]|{\theta} A \ar[d]_f\ar[r]^m|-@{|} & B \ar[d]^g & & & \ar@{}[dr]|{\delta} X \ar[d]_h\ar[r]^p|-@{|} & Y \ar[d]^k \\
        C \ar[r]_n|-@{|} & D & & & C \ar[r]_n|-@{|} & D
    }$$
together with arrows $u$ and $v$ such that $fu = h$ and $gv = k$, there is a unique cell $\gamma\colon p\Rightarrow m$ with source $h$ and target $k$ such that $\theta \gamma = \delta$ holds. A \textbf{restriction} of a niche as at left below is a cartesian cell
    $$\xymatrix{
        A \ar[d]_f & B \ar[d]^g & \ar@{}[dr]|{\leadsto} & & \ar@{}[drr]|{\rho} A \ar[d]_f\ar[rr]^{f_!\otimes n\otimes g^*}|-@{|} && B \ar[d]^g \\
        C \ar[r]_n|-@{|} & D & & & C \ar[rr]_n|-@{|} && D
    }$$
as on the right. Restrictions of this kind are so-called because they provide the restriction, or reindexing, functors for the fibration $\langle\src,\tgt\rangle\colon \D_1\to\D_0\times\D_0$. Thus, dually, an \textbf{opcartesian} cell is an opcartesian arrow for the same functor. An \textbf{extension} of a ``coniche" as on the left below is an opcartesian cell
    $$\xymatrix{
        A \ar[d]_f \ar[r]^m|-@{|} & B \ar[d]^g & \ar@{}[dr]|{\leadsto} & & \ar@{}[drr]|{\xi} A \ar[d]_f\ar[rr]^m|-@{|} && B \ar[d]^g \\
        C & D & & & C \ar[rr]_{f^*\otimes m \otimes g_!}|-@{|} && D
    }$$
as on the right. Put another way, a double category $\D$ is an equipment if, and only if, the functor $\langle \src,\tgt\rangle\colon \D_1\to\D_0\times\D_0$ has all restrictions and extensions for given pairs of ordinary arrows are left adjoint to the corresponding restrictions  \cite[Proposition 9.1.2]{BJ}. Special cases of these restrictions and extensions are of importance. Recall \cite{GP2} that an arrow $f\colon A\to B$ and proarrow $f_!\colon A\hto B$ are \textbf{companions} with unit and counit 
    $$\xymatrix{
        \ar@{}[dr]|{\xi} A \ar[d]_1 \ar[r]^{y_A}|-@{|} & A \ar[d]^f & & & \ar@{}[dr]|{\rho} A \ar[d]_f \ar[r]^{f_!}|-@{|} & B \ar[d]^1 \\
        A \ar[r]_{f_!}|-@{|} & B & & & B \ar[r]_{y_B}|-@{|} & B 
    }$$
if $\xi\otimes \rho=1$ and $\rho\xi=y_f$ both hold. Dually, an arrow $f\colon A\to B$ and proarrow $f^*\colon B\hto A$ are \textbf{conjoint} with unit and counit
    $$\xymatrix{
        \ar@{}[dr]|{\rho} B \ar[d]_1 \ar[r]^{f^*}|-@{|} & A \ar[d]^f & & & \ar@{}[dr]|{\xi} A \ar[d]_f \ar[r]^{y_A}|-@{|} & A \ar[d]^1 \\
        B \ar[r]_{y_B}|-@{|} & B & & & B \ar[r]_{f^*}|-@{|} & A 
    }$$
if $\rho\otimes\xi=1$ and $\rho\xi = y_f$ hold. In the former case $(f,f_!)$ is said to form a companion pair; in the latter case $(f,f^*)$ is a conjoint pair. $f_!$ is a companion of $f$ and $f^*$ is a conjoint of $f$. Companions and conjoints completely describe equipment structure in the sense that $\D$ is an equipment if, and only if, $\D$ has all companions and conjoints. A detailed proof is given in \cite{FramedBicats}. In particular, the manner in which restrictions and extensions can be built from companions and conjoints is suggested in the notation `$f_!\otimes n\otimes g^*$'and `$f^*\otimes m \otimes g_!$'. Restrictions and extensions of external identities will be denoted as `$f_!\otimes g^*$' and `$f^*\otimes g_!$' without the `$y$' to reduce notational clutter.

\example
    $\dblset$, $\prof$, $\Rel(\E)$ are equipments. Companions in $\Rel(\E)$ are given by graphs; conjoints by opgraphs. Extensions in $\Rel(\E)$ are computed by images; restrictions are given by pullback. In particular, any morphism $e\colon E\to B$ in $\E$ is a cover in $\E$ if, and only if, the corresponding cell $y_e$ coming from the external identity is an extension. For the extension 
        $$\xymatrix{
            \ar@{}[dr]|{\xi} A \ar[d]_e \ar[r]^{y_A}|-@{|} & A \ar[d]^e \\
            E \ar[r]_{e^*\otimes e_!}|-@{|} & E
        }$$
    results in a globular cell $\gamma\colon e^*\otimes e_! \Rightarrow y_E$ such that $\gamma\xi = y_e$. But this is computed by an image in $\E$. So, $e$ is a cover if, and only if, the unique globular cell $\gamma$ above is an iso $e^*\otimes e_!\cong y_E$, if and only if, $y_e$ is an extension. 
\endexample

\definition[Cf. \S 4.2, \cite{Schultz}]
    The \textbf{kernel} of a morphism $f\colon A\to B$ is the restriction $\rho$ of the unit on $B$ along $f$. Dually, the \textbf{cokernel} of $f$ is the extension cell $\xi$
        $$\xymatrix{
        \ar@{}[dr]|{\rho} A \ar[d]_f \ar[r]^-{f_!\otimes f^*}|-@{|} & A \ar[d]^f & & & \ar@{}[dr]|{\xi} A \ar[d]_f \ar[r]^-{y_A}|-@{|} & A \ar[d]^f \\
        B \ar[r]_{y_B}|-@{|} & B & & & B \ar[r]_{f^*\otimes f_!}|-@{|} & B
        }$$
    A morphism $e\colon A\to E$ in an equipment is a \textbf{cover} if the canonical globular cell is an iso $e^*\otimes e_!\cong y_E$. Dually, $m\colon E\to B$ is an \textbf{inclusion} if the canonical globular cell is an iso $m_!\otimes m^*\cong y_E$.
\enddefinition

A pithy or sloganish formulation is that a cover is a morphism with trivial cokernel; an inclusion is one with trivial kernel. Notice that these recall the definitions of simple and entire maps in an allegory \cite[\S II,2.13]{FS}. This partly explains the choice of notation for companions and conjoints. Below it will be seen that at least for certain `double categories of relations' $f\mapsto f^*$ is a kind of anti-involution operator similar to the one axiomatized in the definition of an allegory.

Preservation of extensions and restrictions by oplax or lax functors will be important throughout. The main result in this connection is the following.

\proposition \label{prop:Oplax/LaxPreservationResult}
    An oplax double functor preserves extensions; a lax one preserves restrictions. Thus, a pseudo-double functor preserves both.
\endproposition
\proof
    See Proposition 6.4 and its proof in \cite{FramedBicats}.
\endproof

Now, turn to the cartesian structure of relations. Let $\dbl$ denote the 2-category of double categories, pseudo-double functors and (vertical) transformations. The definition and properties of cartesian double categories are given in \S 4 of \cite{aleiferi}.

\definition[Cf. \S 4.2, \cite{aleiferi}] \label{define:CartesianDoubleCategory}
    A double category $\D$ is \textbf{cartesian} if the double functors $\Delta\colon \D\to\D\times \D$ and $\D\to\mathbf 1$ have right adjoints in $\dbl$.
\enddefinition

If $\D$ is cartesian, then $\D_0$ and $\D_1$ both have finite products (Proposition 4.2.2 of \cite{aleiferi}). Additionally $\D$ has ``local products'' in the following sense.

\lemma[Cf. Prop. 4.3.2, \cite{aleiferi}] \label{lemma:CartesianImpliesLocalProducts}
    If $\D$ is a cartesian equipment, then every category $\D(A,B)$ has products.
\endlemma
\proof 
    Given two proarrows $m\colon A\hto B$ and $n\colon A\hto B$, the product in $\D(A,B)$ is given by taking the restriction
        $$\xymatrix{ 
            \ar@{}[dr]|{\rho} A \ar[d]_\Delta \ar@{-->}[r]^-{m\wedge n}|-@{|} & B \ar[d]^{\Delta} \\
            A\times A \ar[r]_{m\times n}|-@{|} & B\times B
        }$$
    along the diagonals. Note that by the construction of restrictions, this is the composite $\Delta_!\otimes (m\times n)\otimes \Delta^*$ which is the formula given in the proof of Theorem 1.6 (ii) of \cite{CW}. The terminal object is the restriction of the identity on $1$ along the unique morphisms from $A$ and $B$ to $1$.
\endproof

Double categories such as $\Span$ and $\prof$ are cartesian equipments. The main example however is relations on a regular category.

\lemma
    $\Rel(\E)$ is cartesian and thus has local products which are computed as intersections of relations.
\endlemma
\proof 
    The underlying category is $\E$, which has a terminal object and products by assumption. The product of relations $R\colon A\hto B$ and $S\colon C\hto D$ is the product
        \[ 
            R\times S\to A\times B\times C\times D\cong A\times C\times B\times D
        \]
    which is again a relation. This defines the required right adjoint to the diagonal double functor. The adjoint for the double functor $\Rel(\E)\to \mathbf 1$ uses the fact that $\E$ has a terminal object. The terminal relation is the evident morphism $1\to 1\times 1$. Since restrictions are computed as pullbacks, the formula in the proof of Lemma \ref{lemma:CartesianImpliesLocalProducts} implies that local products of relations are intersections.
\endproof

The last general property of relations is the interaction between composition, base change and local products. This is the so-called ``Modular Law'' which is part of the definition of an allegory \cite{FS}. It is also closely related to the Frobenius Law. In an allegory, the former is the assertion that there is a valid inequality
    \[
        RT\wedge S \leq (R\wedge ST^\circ)T
    \]
written in the usual compositional order. Here $(-)^\circ$ is the involution coming with the allegory structure which gives the right adjoints for maps. Since $(-)^\circ$ satisfies $((-)^\circ)^\circ = 1$, there is the special case 
    \[
        RT\wedge S \leq (R\wedge ST)T^\circ.
    \]
Insofar as an opgraph is a right adjoint to a graph in a double category $\Rel(\E)$, it might be expected that there is some such relationship governing their interaction with local products. And indeed this is the case. The relationship is the so-called \textbf{modular law}, formulable for any locally posetal cartesian equipment, stating that
    \begin{equation} \label{equation:ModularLaw}
        f^*\otimes R \wedge S \leq f^*\otimes (R\wedge f_!\otimes S) 
    \end{equation}
holds for any morphism $f\colon A\to B$ and proarrows $R\colon A\hto X$ and $S\colon B\hto X$. This statement looks ahead to the further developments of \S \ref{section:CharacterizationTheorem} where it will be seen in Proposition \ref{proposition:InvolutionOperatorExists} that any double category of relations has an involution $(-)^\circ$ with in particular $(f_!)^\circ = f^*$.

\lemma \label{prop:RelSatisfiesFrobenius}
    Local products in $\Rel(\E)$ satisfy the Modular Law.
\endlemma
\proof 
    It suffices to look at the bicategorical fragment of $\Rel(\E)$ and prove the modular law in this context since all the relevant structure is the same. A complete proof for this case is that of Proposition A3.1.5 in \cite{Elephant}.
\endproof

\remark
    In fact the inequality is an equality using a dual formulation of the law. This is equivalent to the \textbf{Frobenius Law} for the hyperdoctrine \cite{hyper} given by
        \[
            \Rel(\E)(-,-)\colon (\D_0\times \D_0)^{op} \to \cat \qquad (A,B) \mapsto \Rel(\E)(A,B)
        \]
    stating that extension as a left adjoint to restriction partially distributes over local products. 
\endremark

\section{`Double Categories of Relations'}
\label{section:DoubleCategoriesOfRelations}

The connection to cartesian bicategories can now be made explicit. This leads to a definition of `double categories of relations' below. The single quotations follow the conventions of \cite{CW} to distinguish the axiomatic notion presented below from double categories defined to be $\Rel(\E)$ for some regular category $\E$.

Recall more precisely that a \textbf{cartesian bicategory} as in \cite{CW} is a locally posetal bicategory with a monoidal tensor for which every object is a commutative comonoid, every morphism $X\to Y$ is a lax comonoid homomorphism, and every comonoid structure map and counit has a right adjoint. It is further required that this comonoid structure on a given object is the unique one having these right adjoints. Recall \cite{GP} that any double category $\D$ has a \textbf{horizontal bicategory} $\mathcal H(\D)$ formed by restricting to the globular cells, that is, the cells whose external source and target morphisms are identities.

\proposition \label{prop:HorizontalBicatIsCartesianBicat}
    The horizontal bicategory of a locally posetal cartesian equipment is a cartesian bicategory.
\endproposition
\proof
    First note that the pseudo-double functor $\times\colon \D\times \D\to \D$ induces a suitable homomorphism of horizontal bicategories. Now, fix an object $B$. The comultiplication and counit on $B$ are defined to be the proarrows arising in the extension cells
        $$\xymatrix{
            \ar@{}[dr]|{\xi} B \ar[d] \ar[r]^y|-@{|} & B \ar[d] & & \ar@{}[dr]|{\xi}B \ar[d] \ar[r]^y|-@{|} & B \ar[d] \\
            B \ar[r]_-{\Delta_!}|-@{|} & B\times B & & B\ar[r]_{\tau_!}|-@{|} & 1.
        }$$
    Note that these could equally well be given by restriction cells. In any case, both have right adjoints $\Delta^*$ and $\tau^*$ given by the equipment structure. The double functor $B\times -\colon \D\to\D$ preserves extensions by Proposition \ref{prop:Oplax/LaxPreservationResult}. So, the two composites
        $$\xymatrix{
            \ar@{}[dr]|{\xi} B \ar[d] \ar[r]^-{\Delta_!}|-@{|} & \ar@{}[dr]|{\xi} B\times B \ar[d] \ar[r]^y|-@{|} & B\times B \ar[d]^{1\times \Delta} & & \ar@{}[dr]|{\xi} B \ar[d] \ar[r]^-{\Delta_!}|-@{|} & \ar@{}[dr]|{\xi} B\times B \ar[d] \ar[r]^y|-@{|} & B\times B\ar[d]^{\Delta\times 1} \\
            B\ar[r]_-{\Delta_!}|-@{|} & B\times B \ar[r]_-{y\times\Delta_!}|-@{|} & B\times (B\times B) & & B\ar[r]_-{\Delta_!}|-@{|} & B\times B \ar[r]_-{\Delta_!\times y}|-@{|} & (B\times B)\times B
        }$$
    compute the same extension, meaning that the comultiplication law $(y\times \Delta_!)\Delta_! = (\Delta_!\times y)\Delta_!$ must hold. Similarly for the unit law. It remains to see that each morphism $p\colon B\hto C$ of $\mathcal H(\D)$ is a lax comonoid homomorphism. But this is straightforward by the construction of the comultiplication and counit morphisms. For example, given the two composite cells
        $$\xymatrix{
            \ar@{}[dr]|{\xi} B \ar[d] \ar[r]^y|-@{|} & \ar@{}[dr]|{\delta} B \ar[d] \ar[r]^p|-@{|} & C \ar[d]^\Delta & & \ar@{}[dr]|{1} B \ar[d] \ar[r]^-p|-@{|} & \ar@{}[dr]|{\xi} C \ar[d] \ar[r]^y|-@{|} & C \ar[d]^\Delta \\
            B\ar[r]_{\Delta_!}|-@{|} & B\times B \ar[r]_-{p\times p}|-@{|} &C\times C & & B\ar[r]_p|-@{|} & C \ar[r]_-{\Delta_!}|-@{|} &C\times C 
        }$$
    the rightmost is an extension, meaning that
        \begin{equation} \label{equation:LaxComonoidHomomorphism}
            p\otimes\Delta_!\leq \Delta_!\otimes (p\times p)  
        \end{equation}
    holds as required. The identity rule is similar. Finally, this comonoid structure needs to be seen to be unique. The argument is a hybrid of those proving Lemma A3.2.3 in \cite{Elephant} and Corollary 2.1.6 in \cite{aleiferi}. Suppose that $(B,d_B,e_B)$ is another comonoid structure on $B$ for which $d$ and $e$ have right adjoints $d^*$ and $e^*$. Since $y_1$ is terminal in $\D_1$, there is a unique cell 
        $$\xymatrix{
            \ar@{}[dr]|{\exists\,!} B\ar[d] \ar[r]^e|-@{|} & 1 \ar[d] \\
            1 \ar[r]|-@{|} & 1
        }$$       
    meaning that $e\leq \tau_!$ holds since $\tau_!$ is equivalently given by a cartesian cell over $y_1$. Using the unit and counit of the adjunctions, it follows that 
        \[
            \tau_!\leq \tau_!e^*e^*\leq \tau_!\tau^*e\leq e  
        \]
    proving that $e=\tau_!$ holds. Now, use the argument of Corollary 2.1.6 in \cite{aleiferi} which shows that $d$ and $\Delta_!$ are coequalized by the projections from $B\times B$ and are thus equal.
\endproof

What makes a cartesian bicategory a `bicategory of relations' is a further discreteness condition, namely, that for each object $B\in \mathbb B$, the corresponding comonoid structure morphism $\Delta$ and its right adjoint $\Delta^*$ satisfy the equation
    \[
        \Delta^*\otimes \Delta = (1\times \Delta)\otimes (\Delta^*\times 1).  
    \]
This is called the \textbf{Frobenius Law} in \cite{CW}. It now makes sense to define a `double category of relations' as a locally posetal cartesian equipment satisfying the analogous discreteness condition. 

\definition \label{define:DoubleCategoryOfRelations}
    A `\textbf{double category of relations}' is a locally posetal cartesian equipment in which every object is \textbf{discrete} in the sense that the Frobenius Law
        \begin{equation}
            \Delta^*\otimes \Delta_! = (1\times \Delta_!)\otimes (\Delta^*\times 1)\qquad 
        \end{equation}
    holds. The single quotations in `double category of relations' will distinguish this axiomatic notion from double categories equal to $\Rel(\E)$ for some $\E$.
\enddefinition

Consequently, the horizontal bicategory of any `double category of relations' is indeed a `bicategory of relations'. It could be viewed as a generalized compact closed category \cite{KellyLaplaza}. See Theorem 2.4 in \cite{CW} for more on that point. For the present purposes, a number of consequences follow from this connection with `bicategories of relations'. For example, the Modular Law \ref{equation:ModularLaw} follows from the `double category of relations' axioms. This will be used in the proof of Theorem \ref{theorem:D0isRegular} below.

\corollary
    In any `double category of relations', the modular law \ref{equation:ModularLaw} holds.
\endcorollary
\proof
    Remark 2.9(ii) of \cite{CW} shows that any bicategory of relations satisfies the modular law. Thus, the modular law certainly holds for the horizontal bicategory $\mathcal H(\D)$ of any `double category of relations'. However, since these adjoints and products are inherited from $\D$, this means that the modular law holds in $\D$ too.
\endproof

It can also be seen that local products are idempotent. This follows by first establishing that diagonal morphisms are inclusions. Note that this only requires the cartesian equipment structure.

\corollary \label{prop:DiagonalsAreInclusionsIESeparability}
    In a locally posetal cartesian equipment, the diagonal morphisms $A\to A\times A$ are inclusions.
\endcorollary
\proof
    Fix an object $A$. By Proposition \ref{prop:HorizontalBicatIsCartesianBicat} this has the structure of a commutative comonoid in the horizontal bicategory. The goal is to prove that $\Delta_!\otimes \Delta^* = y$ via the canonical map. If $\epsilon_!\colon A\hto 1$ is the counit of $A$ viewed as a comonoid, the unit $1\leq \epsilon_!\otimes \epsilon^*$ yields an inequality
        \begin{equation}
            \Delta_!\otimes \Delta^* \leq \Delta_!\otimes (y\times \epsilon_!)\otimes (y\times \epsilon^*)\otimes \Delta^*
        \end{equation}
    However, the right side is $y$. For the composite
        $$\xymatrix{
            \ar@{}[dr]|{\rho} A \ar[d]_\Delta \ar[r]^-{\Delta_!}|-@{|} & \ar@{}[dr]|{1} A\times A \ar[d] \ar[r]^{y\times\epsilon_!}|-@{|} &A\times 1 \ar[d] \\
            \ar@{}[drr]|{1\times \rho} A\times A \ar[d] \ar[r]_y|-@{|} &A\times A \ar[r]_{y\times\epsilon_!}|-@{|} &A\times 1 \ar[d] \\
            A\times 1 \ar[rr]_y|-@{|} && A\times 1
        }$$
    is a restriction cell, implying that $\Delta_!\otimes (y\times \epsilon_!)\cong y$. Similarly for $(y\times\epsilon^*)\otimes \Delta^*$. Since $y\leq \Delta_!\otimes \Delta^*$ always holds (the right side is a restriction), this proves the result.
\endproof

\corollary \label{lemma: local products idempotent}
    In any `double category of relations', local products of companions $p_!$ are idempotent in the sense that $p_! = p_!\wedge p_!$ holds canonically. Equivalently, each diagonal cell
        $$\xymatrix{
            \ar@{}[dr]|{\Delta_{p_!}} A \ar[d]_{\Delta_A} \ar[r]^{p_!}|-@{|} & B \ar[d]^{\Delta_B} \\
            A\times A  \ar[r]_{p_!\times p_!} & B\times B
        }$$
    is cartesian. Dually, local products of the right adjoints $p^*$ are idempotent.
\endcorollary
\proof 
    There are canonical isomorphisms
        \begin{equation} \label{equation:LaxComonoidHomomIsStrict}
            p_!\otimes \Delta_! \cong \Delta_!\otimes (p\times p)_! \cong \Delta_!\otimes (p_!\times p_!)
        \end{equation}
    The leftmost iso is from the fact that two sides compute the same restriction, namely, that of the morphism $\Delta_Bp = (p\times p)\Delta_A$. The rightmost is from the fact that the product functor $\times \colon \D\times\D\to\D$ is pseudo, hence preserves restrictions by Proposition \ref{prop:Oplax/LaxPreservationResult}. As a result, the isomorphism
        $$\xymatrix{ 
            \ar@{}[dr]|{1} A \ar@{=}[d] \ar[r]^{p_!}|-@{|} & \ar@{}[drr]|{\cong} B \ar@{=}[d] \ar[rr]^y|-@{|} & &  B \ar@{=}[d] \\
            \ar@{}[drr]|{\cong} A \ar@{=}[d] \ar[r]_{p_!}|-@{|} & B \ar[r]^-{\Delta_!}|-@{|} & \ar@{}[dr]|{1} B\times B \ar@{=}[d] \ar[r]^-{\Delta^*}|-@{|} & B \ar@{=}[d] \\
            A \ar[r]_-{\Delta_!}|-@{|} & A\times A \ar[r]_{p_!\times p_!}|-@{|} & B\times B \ar[r]_-{\Delta^*}|-@{|} & B
        }$$
    establishes the result by the construction of local products from the proof of Lemma \ref{lemma:CartesianImpliesLocalProducts}. The topmost iso is the fact that $\Delta_B$ is an inclusion by Corollary \ref{prop:DiagonalsAreInclusionsIESeparability}. The dual case is analogous. 
\endproof

\remark
    Equation \ref{equation:LaxComonoidHomomIsStrict} is saying that each $p_!$ is a comonoid homomorphism. Since $p_!$ has a right adjoint $p^*$ this is basically the characterization of so-called ``maps" in Lemma 2.5 of \cite{CW}. 
\endremark

\section{Tabulators and Functional Completeness}
\label{section:RelsTabulatorsAreStrongandMonic}

Allegories and cartesian bicategories have their respective notions of a ``tabulation'' of a given arrow. In the context of double categories, the tabulator of a proarrow is a kind of finite limit.

\definition[\cite{GP}]
    A double category $\D$ has \textbf{tabulators} if $y\colon \D_0\to \D_1$ has a right adjoint $\top\colon \D_1\to \D_0$ in $\dbl$. The \textbf{tabulator} of a proarrow $m\colon A\hto B$ is the object $\top m$ together with the counit cell $\top m\Rightarrow m$. Denote the external source and target by $l$ and $r$.
\enddefinition

\lemma \label{lemma:RelHasTabs}
    $\Rel(\E)$ has tabulators. The unit of $y\dashv \top$ is iso. Equivalently, $y$ is fully-faithful.
\endlemma
\proof 
    Define $\top\colon \Rel(\E)_1\to\Rel(\E)_0$ by sending $R\to A\times B$ to $R$ with the evident assignment on arrows. In other words, $\top$ takes the apex of spans and morphisms between them. The component of the counit at $R\to A\times B$ is the cell given by the morphism of relations
        $$\xymatrix{
            R \ar[d]_1\ar[r]^-{\Delta} & R\times R \ar[d]^{d\times c} \\
            R \ar[r]_-{\langle d,c\rangle} & A\times B
        }$$
    On the other hand, the unit is up to iso the identity map on a given object $A$. That is, $y$ takes the diagonal $A\to A\times A$ and then $\top$ takes the apex $A$, meaning that $1\cong \top y$ canonically. By a general result \cite[IV.3.1]{MacLane} this is equivalent to the statement that $y$ is fully faithful.
\endproof

\definition[\S 4.3.7, \cite{aleiferi}]
    A double category $\D$ is \textbf{unit-pure} if the external identity $y\colon \D_0\to \D_1$ is fully faithful. 
\enddefinition

\example
    $\dblset$, $\Span(\E)$ and $\Rel(\E)$ are all unit-pure whereas $\prof$ is not.
\endexample

Since $y$ is always faithful, technically all that is required in the definition is ``full." Tabulators in $\Rel(\E)$ are additionally ``strong" in the following sense since extensions are given by taking image factorizations.

\definition
    The tabulator $\langle l,r\rangle\colon \top m \to A\times B$ of a proarrow $m\colon A\hto B$ is \textbf{strong} if $m$ is the cokernel of its tabulator in the sense that $m\cong l^*\otimes r_!$ holds canonically.
\enddefinition

In a unit-pure equipment with strong tabulators, inclusions are precisely the monic arrows.

\lemma \label{lemma:PropertiesOfInclusionsUnitPure}
    If $y\colon \D_0\to \D_1$ is fully faithful, then 
        \begin{enumerate}
            \item inclusions are monic; 
            \item if tabulators are strong, then monic arrows are inclusions.
        \end{enumerate}
\endlemma
\proof 
    (1) Assume that $f\colon A\to B$ is an inclusion and take arrows $u,v\colon X\rightrightarrows A$ such that $fu=fv$. Using the fact that $f$ is an inclusion and that $y_{fu} = y_{fv}$, it follows that $y_u=y_v$ holds. From this $u=v$ since $y$ is faithful. 

    (2) On the other hand, given a monic $f\colon A\to B$, take the tabulator of its kernel. Denote the legs by $l$ and $r$. Since $B$ is the tabulator of $y_B$, it follows that $fl=fr$ must hold. But then $l=r$ follows since $f$ is monic. Now, since tabulators are strong, the kernel of $f$ is isomorphic to the cokernel of its tabulator. Thus, there are unique (globular) cells
        \[
          y_A \Rightarrow f_!\otimes f^* \cong l^*\otimes l_!\Rightarrow y_A.
        \] 
    Since $y$ is fully faithful, these must compose to the identity. The other composite is also the identity by uniqueness of lifts. Thus, $f_!\otimes f^*\cong y_A$ holds proving that $f$ is an inclusion.
\endproof 

Tabulators in $\Rel(\E)$ are then monic in the following sense.

\proposition \label{proposition:TabulatorsInRelAreMonic}
    The legs of the tabulator of a proarrow $R\colon A\hto B$ in $\Rel(\E)$ are jointly monic and the iso
        \begin{equation} \label{equation:MonicTabulators}
            l_!\otimes l^*\wedge r_!\otimes r^*\cong y
        \end{equation}
    holds canonically.
\endproposition
\proof
    The proof of Lemma \ref{lemma:RelHasTabs} shows that any relation is its own tabulator. Thus, the first statement is trivial owing to the fact that a proarrow is just a monic arrow $R\to A\times B$. For the second statement, the restriction gives the kernel of the relation:
        $$\xymatrix{
            \ar@{}[dr]|{\rho} R \ar[d]_{\langle l,r\rangle} \ar[r]^{\langle l, r\rangle_!\otimes \langle l,r\rangle^*}|-@{|} & R \ar[d]^{\langle l,r\rangle} \\
            A\times B \ar[r]_y|-@{|} & A\times B
        }$$
    But the composite
        $$\xymatrix{ 
            \ar@{}[drr]|{\rho} R \ar[d]_\Delta \ar[rr]^{l_!\otimes l^*\wedge r_!\otimes r^*}|-@{|} && R \ar[d]^\Delta \\
            \ar@{}[drr]|{\rho} R \times R \ar[d]_{l\times r} \ar[rr]_{l_!\otimes l^*\times r_!\otimes r^*}|-@{|} && R\times R \ar[d]^{l\times r} \\
            A\times B \ar[rr]_{y\times y}|-@{|} && A\times B 
        }$$
    computes the same restriction. So, if one is $y$, then the other is too and conversely.
\endproof

\definition[Cf. \S 3, \cite{CW}] \label{define:FunctionallyComplete}
    A cartesian equipment is \textbf{functionally complete} if it has tabulators and these are strong and monic in the sense that if $R\colon A\hto B$ is any proarrow with tabulator $\top R\to A\times B$ with legs $l$ and $r$, then the equations
        \[
             R\cong l^*\otimes r_! \qquad\text{and}\qquad  l_!\otimes l^*\wedge r_!\otimes r^*\cong y
        \]
    both hold.
\enddefinition

Consequently, in any unit-pure functionally complete cartesian equipment, the legs of any tabulator are genuinely jointly monic by Lemma \ref{lemma:PropertiesOfInclusionsUnitPure}.

\section{Characterizing Spans}
\label{section:SpansReview}

The next few sections are devoted to the first part of the main result. This is an answer to the question of the conditions under which a double category $\D$ is equivalent to $\Rel(\D_0)$ where $\D_0$ is a regular category in the sense of Definition \ref{define:RegularCat}. The starting point is the now well-established treatment of spans in a cartesian category. This has its roots in a result of \cite{Niefield}. Namely, a double category $\D$ admits a normalized oplax/lax adjunction to $\Span(\D_0)$ if, and only if, it is an equipment with tabulators.

\theorem[Theorems 5.5/5.6, \cite{Niefield}] \label{theorem:NiefieldTheorem}
    Let $\D$ denote a double category with pullbacks. The following are equivalent:
        \begin{enumerate}
            \item There is an oplax/lax adjunction $F\colon \Span(\D_0) \rightleftarrows \D \colon G$ where $F$ is normal and equal to the identity on $\D_0$.
            \item $\D$ has all companions, conjoints and tabulators.
        \end{enumerate}
\endtheorem

\proof
    That $\D$ has tabulators and is an equipment allows construction of the oplax and lax functors. $G$ is defined on proarrows by taking a tabulator; $F$ is defined on proarrows by taking an extension of a coniche given by a span. 
\endproof

The main result of \cite{aleiferi} gives equivalent conditions under which such a normalized oplax/lax adjunction is a strong equivalence. These extra conditions are that $\D$ is unit-pure, cartesian and possesses certain internal Eilenberg-Moore objects defined in the following way.

\definition[Cf. \S 5.3 of \cite{aleiferi}]
    A \textbf{copoint} of an proarrow $m\colon A\hto A$ in $\D$ is a cell 
        $$\xymatrix{ 
            \ar@{}[dr]|{\gamma} A\ar@{=}[d] \ar[r]^m|-@{|} & A \ar@{=}[d] \\
            A \ar[r]_{y_A}|-@{|} & A.
        }$$
    Let $\copt(\D)$ denote the category of pairs $(m,\gamma)$ where $\gamma$ is a copoint of the endoproarrow $m$. The morphisms $(m,\gamma) \to (n,\epsilon)$ are cells $\theta\colon m\Rightarrow n$ of $\D$ such that $\epsilon\theta = \gamma$ holds. A double category $\D$ \textbf{admits Eilenberg-Moore objects for copointed endomorphisms} if the inclusion $\D_0 \to\copt(\D)$ has a right adjoint.
\enddefinition

The characterization of spans is then the following. Its proof is the topic of \S 5 of the reference and so will not be reproduced here. The Beck-Chevalley condition appearing in the third equivalent condition has not yet been discussed but is stated for an equipment in Definition \ref{define:BeckChevalley} below. Its discussion is postponed only because it is not part of the template for the present results.

\theorem[Theorem 5.3.2, \cite{aleiferi}]
    For a double category $\D$ the following are equivalent:
    \begin{enumerate}
        \item $\D$ is equivalent to $\Span(\E)$ for some finitely-complete category $\E$.
        \item $\D$ is a unit-pure cartesian equipment admitting Eilenberg-Moore objects for copointed endoproarrows.
        \item $\D_0$ has pullbacks satisfying the strong Beck-Chevalley condition and the canonical functor
            \[ 
                \Span(\D_0) \to \D
            \]
        is an equivalence of double categories.
    \end{enumerate}
\endtheorem
\proof
    This is stated and proved completely in \S 5.3 of the reference.
\endproof

The development for relations will follow this pattern. Namely, start with conditions equivalent to the existence of an oplax/lax adjunction and isolate the further conditions under which such an adjunction is a strong equivalence of double categories.

\section{Characterizing Relations: Conditions for Oplax/Lax Adjunction}
\label{section:OplaxLaxAdjunctionConditions}

First develop the relation version of Niefield's Theorem \ref{theorem:NiefieldTheorem} quoted above. This appears as Theorem \ref{theorem:ExistenceOplaxLaxAdjunction} below. Consider first some necessary conditions.

\lemma \label{lemma:NecessaryConditionsForOplaxLaxAdjunction}
    Let $\D$ denote a double category where $\D_0$ is a regular category. Suppose that $F\colon \Rel(\D_0)\rightleftarrows \D\colon G$ is a normalized oplax/lax adjunction that is the identity on $\D_0$. It then follows that $\D$
        \begin{enumerate}
            \item is an equipment;
            \item has monic tabulators;
            \item the unit $1\Rightarrow\top y$ is an iso;
            \item and for each cover $e\colon A\twoheadrightarrow  E$, the cell
                $$\xymatrix{ \ar@{}[dr]|{y_e} A \ar@{->>}[d]_e \ar[r]^{y_A}|-@{|} & A \ar@{->>}[d]^{e} \\
                E \ar[r]_{y_E}|-@{|} & E
                }$$
            is an extension in $\D$.
        \end{enumerate}
\endlemma
\proof 
    Take $f\colon A\to B$ in $\D_0$. The graph and opgraph give the companion and conjoint in $\Rel(\D_0)$. The images under $F$ give the corresponding companion and conjoint in $\D$ making it an equipment. Since oplax functors preserve extensions and every cell
        $$\xymatrix{
            \ar@{}[dr]|{e} A \ar@{->>}[d]_e \ar[r]^\Delta|-@{|} & A \ar@{->>}[d]^e \\
            E \ar[r]_\Delta|-@{|} & E
        }$$
    is one in $\Rel(\D_0)$, the corresponding image under $F$ is an extension, making $y_e$ an extension in $\D$ since it is isomorphic to $Fe$ by normalization. Existence of tabulators results from the fact that the composite 
        \[ 
            \D_1 \xrightarrow[]{G_1}\Rel(\D_0)_1 \xrightarrow[]{apex} \D_0
        \]
    is a right adjoint for $y\colon \D_0\to \D_1$. By normalization of $G$, its unit is an isomorphism. 
\endproof

Now, it can be seen that these conditions are also sufficient. It is worth doing the details of these constructions since they show precisely what is required for the resulting adjunction to be a strong adjoint equivalence. Assume throughout that $\D_0$ is regular.

\lemma \label{lemma:SufficiencyForOplaxF}
    If $\D$ has companions and conjoints and $y_e$ is an extension for each cover $e$, then the identity $1\colon \D_0\to\D_0$ extends to an opnormal oplax functor $F\colon\Rel(\D_0) \to \D$.
\endlemma
\proof 
    For an relation $R\to A\times B$, take the image $FR$ in $\D$ to be the proarrow $A\hto B$ arising in the canonical extension
        $$\xymatrix{
            \ar@{}[dr]|{\xi_R} R \ar[d]_d \ar[r]^y|-@{|} & R \ar[d]^c \\
            A \ar[r]_{d^*\otimes c_!}|-@{|} & B
        }$$
    That the cell is opcartesian gives the arrow assignment, hence by uniqueness properties a functor $F_1\colon \Rel(\D_0)_1\to\D_1$. Comparison cells for composition are given using the extension property of the composite cell $\xi y_e$. That is, they arise as in the lower-left corner of the diagram:
        $$\xymatrix{
            R\times_BS \ar@{=}[d] \ar[rrrr]^{R\times_BS}|-@{|} &  & &  & R\times_BS \ar@{=}[d] & & \ar@{}[drr]|{y_e} R\times_BS \ar@{->>}[d]_e \ar[rr]^{R\times_BS}|-@{|} & & R\times_BS \ar@{->>}[d]^e \\
            \ar@{}[dr]|{y_p} R\times_BS \ar[d]_p \ar[rr]^{R\times_BS}|-@{|} & & R\times_BS \ar[dl]_p \ar[dr]^q \ar[rr]^{R\times_BS}|-@{|} &\ar@{}[dr]|{y_q} & R\times_BS \ar[d]^q \ar@{}[drr]|{=}& & \ar@{}[drr]|{\xi} R\otimes S \ar[d]_d \ar[rr]_{R\otimes S}|-@{|} & & R\otimes S \ar[d]^c \\
            \ar@{}[dr]|{\xi_R} R \ar[d]_d \ar[r]^R|-@{|}& R \ar[d]^c & & \ar@{}[dr]|{\xi_S} S \ar[d]_d \ar[r]_{S}|-@{|}& S \ar[d]^c &  &   \ar@{}[drr]|{\exists !\,\phi_{R,S}} A\ar@{=}[d] \ar[rr]_{d^*\otimes c_!}|-@{|}& & C \ar@{=}[d]   \\
            A \ar[r]_{d^*\otimes_R c_!}|-@{|} & B \ar@{=}[rr] & & B \ar[r]_{d^*\otimes_S c_!}|-@{|} & C & & A \ar[r]_{d^*\otimes_R c_!}|-@{|}& B \ar[r]_{d^*\otimes_S c_!}|-@{|} & C 
        }$$
    In general these are not invertible. The coherence laws for an oplax functor follow by the fact that all the cells are defined using the uniqueness clause of the lifting property of opcartesian cells.
\endproof

\lemma \label{lemma:SufficiencyForLaxG}
    If $\D$ has tabulators whose legs are jointly monic, then $1\colon \D_0\to\D_0$ extends to a lax functor $G\colon \D\to\Rel(\D_0)$. It is normal if, and only if, $y$ is fully faithful.
\endlemma
\proof 
    Write $\top\colon \D_1\to\D_0$ for the right adjoint to $y\colon \D_0\to\D_1$. For a proarrow $p\colon A\hto B$ in $\D$, define the image $Gp$ in $\Rel(\D_0)$ to be the inclusion $Tp\to A\times B$ given by the tabulator of $p$. By the universal property of tabulators this induces a functorial arrow assignment yielding the required functor $G_1\colon \D_1\to \Rel(\D_0)$. Externally this is lax-functorial. For composable proarrows $p\colon A\hto B$ and $q\colon B\hto C$, the induced morphism
        $$\xymatrix{
            \top p\otimes \top q \ar@{-->}[d]_{\gamma_{p,q}} \ar[r] & A\times C \ar@[=][d] \\
            \top (p\otimes q) \ar[r] & A\times C,
        }$$
    which exists by orthogonality of the factorization system on the regular category $\D_0$, gives the required laxity cell $\gamma\colon Tp\otimes Tq\Rightarrow T(p\otimes q)$. The naturality and associativity conditions for a lax functor follow by the uniqueness of image factorizations, the fact that tabulators are jointly inclusions, and the fact that regular epis are pullback-stable. Unit comparison cells are induced again from the universal property of tabulators; given an object $A$ the induced morphism
        $$\xymatrix{
            A \ar@{-->}[d]_{\gamma_A} \ar[r]^-\Delta & A\times A \ar[d]^{1\times 1} \\
            \top y_A \ar[r] & A\times A
        }$$
    defines the required cell $\gamma_A\colon y_A\Rightarrow \top y_A$. Note that $G$ is normalized if, and only if, $y\colon \D_1\to\D_0$ is fully faithful.   
\endproof

\proposition
    If $\D$ is an equipment with tabulators whose legs are jointly monic and where $y_e$ is an extension for each cover $e$, then the functors $F_1\colon \Rel(\D_0)_1\rightleftarrows \D_1\colon G_1$ of the previous lemmas form an adjunction $F_1\dashv G_1$.
\endproposition
\proof 
    Develop the unit $\eta\colon 1\Rightarrow G_1F_1$. Starting with a relation $R\rightarrowtail A\times B$, take the canonical extension and then its tabulator. By the universal property of the tabulator, there is a unique morphism $R\to \top (d^*\otimes_Rc_!)$ fitting into
        $$\xymatrix{
            R \ar[d]_{\eta_R} \ar[r] & A\times B \ar@{=}[d]\\
            \top (d^*\otimes_Rc_!) \ar[r] & A\times B
        }$$
    making a morphism of relations. Take this to be the component $\eta_R$. These are natural in $R$ by the uniqueness aspect of the universal property of tabulators. On the other hand, components of the counit $\epsilon\colon F_1G_1\Rightarrow 1$ are given in the following way. For a given proarrow $p\colon A\hto B$, the proarrow $F_1G_1p$ is the extension of the image of the tabulator $\top p$. The counit component $\epsilon_p$ arises as in the the left bottom corner of the diagram
        $$\xymatrix{ 
            \ar@{}[dr]|{\xi_{\top p}} \top p \ar[d]_{l} \ar[r]^y|-@{|} & \top p \ar[d]^r &  &  & \ar@{}[ddr]|{\tau_p} \top p \ar[dd]_l \ar[r]^y|-@{|} &  \top p \ar[dd]^r\\
            \ar@{}[dr]|{\exists\,!} A \ar@{=}[d] \ar[r]_{d^*\otimes_{\top p}c_!}|-@{|} & B \ar@{=}[d] & \ar@{}[r]|{=} & & & \\
            A \ar[r]_p|-@{|} & B &  & & A \ar[r]_p|-@{|} & B 
        }$$
    since the extension $\xi_{\top p}$ is opcartesian. Again this is natural in $p$ by functoriality of tabulators and uniqueness clauses of universal properties. Triangle identities follow by construction. For example, given a relation $R\rightarrowtail A\times B$, verify that $\epsilon_{F_1R}F_1\eta_R = 1$ holds. There are equalities of cells
        $$\xymatrix{
            \ar@{}[dr]|{\xi_R} R \ar[d]_d \ar[r]^y|-@{|} & R \ar[d]^c & & & \ar@{}[dr]|{y_{\eta_R}} R \ar[d]_{\eta_R} \ar[r]^y|-@{|} & R \ar[d]^{\eta_R} & & & \ar@{}[dddr]|{\xi_R} R \ar[ddd]_d \ar[r]^y|-@{|} & R \ar[ddd]^c \\
            \ar@{}[dr]|{F_1\eta_R} A \ar@{=}[d] \ar[r]|-@{|} & B \ar@{=}[d]  & \ar@{}[dr]|{=} & & \ar@{}[dr]|{\xi_{\top(d^*\otimes c_!)}} \top(d^*\otimes c_!) \ar[d]_d \ar[r]|-@{|} & \top(d^*\otimes c_!) \ar[d]^c & \ar@{}[dr]|{=} & & & \\
            \ar@{}[dr]|{\epsilon} A \ar@{=}[d] \ar[r]|-@{|} & B \ar@{=}[d] & & & \ar@{}[dr]|{\epsilon} A \ar@{=}[d] \ar[r]|-@{|} & B \ar@{=}[d] & & & & \\
            A \ar[r]_{d^*\otimes_R c_!}|-@{|} & B & & & A \ar[r]_{d^*\otimes_R c_!}|-@{|} & B & & & A \ar[r]_{d^*\otimes_R c_!}|-@{|} & B
        }$$
    by construction. The leftmost holds by the definition of $F_1\eta_R = d_!\eta_Rc^*$; the right holds by construction of $\epsilon_{d_!Rc^*}$. Now, the composite in the lower left is $\epsilon_{F_1R}F_1\eta_R$. It must be an identity since $\xi_R$ occurring on both sides is an extension. Verifying the other triangle identity is a similar kind of argument but more straightforward.    
\endproof

\theorem \label{theorem:ExistenceOplaxLaxAdjunction}
    For a double category $\D$ where $\D_0$ is regular, the identity $\D_0\to \D_0$ extends to an oplax/lax adjunction $F\colon \Rel(\D_0)\rightleftarrows \D\colon G$ if, and only if, 
        \begin{enumerate}
            \item $\D$ is a unit-pure equipment;
            \item has monic tabulators; 
            \item $y_e$ is an extension for each cover $e$.
        \end{enumerate}
\endtheorem 
\proof 
    There remains only to verify the remaining conditions of an oplax/lax adjunction. These are those of (d) in \S 3.2 of \cite{GP2}. Given composable relations $R\to A\times B$ and $S\to B\times C$, the components of $\eta$ need to be coherent with external composition and laxity cells. But the morphisms of relations on each side of
        $$\xymatrix{ 
            R\otimes S \ar[d]_{\eta_{R\otimes S}}\ar[r]& A\times C\ar@{=}[d] & & R\otimes S \ar[d]_{\eta_R\otimes\eta_S}\ar[r] & A\times C \ar@{=}[d] \ar@{=}[d] \\
            \top(d^*\otimes_{R\otimes S}c_!) \ar[d]_{\top\phi} \ar[r] & A\times C \ar@{=}[d] & & \top (d^*Rc_!)\otimes \top (d^*Sc_!) \ar[d]_{\gamma} \ar[r] & A\times C \ar@{=}[d] \\
            \top(d^*\otimes_{R\otimes S}c_!) \ar[r] & A\times C & & \top(d^*\otimes_{R\otimes S}c_!) \ar[r] & A\times C
        }$$
    are the same by the uniqueness of image factorizations. Similarly, that components of $\epsilon$ are coherent with external composition follows by the construction of $\phi$ and the uniqueness property of cells induced by opcartesian cells.
\endproof

\section{Characterizing Relations: Conditions for Adjoint Equivalence}
\label{section:ExistenceOplaxLaxAdjEquiv}

The constructions from the previous subsection lead to the main result characterizing the existence of an adjoint equivalence $\D\simeq \Rel(\E)$. The proofs of the preliminaries to Theorem \ref{theorem:ExistenceOplaxLaxAdjunction} and some extra streamlining reveal two further conditions guaranteeing an adjoint equivalence. Namely, these are that tabulators are strong and that every relation is a tabulator of its cokernel.

To start, recall that a ``strong" adjoint equivalence of double categories is an oplax/lax adjoint equivalence where both functors are pseudo. In the present development, the proof above shows that this amounts to a Beck-Chevalley condition and the requirement that tabulators are in a particular sense ``functorial." In logic, Beck-Chevalley is the condition, roughly speaking, that substitution commutes with quantification \cite[\S 1.8]{BJ}. Categorically, this is to ask that certain adjoints partially commute. Double categorically this is expressed by the following. 

\definition[\S 13, \cite{FramedBicats}; \S 5.2, \cite{aleiferi}] \label{define:BeckChevalley}
    An equipment $\D$ satisfies the \textbf{Beck-Chevalley condition} if for any pullback square
        $$\xymatrix{ \cdot \ar[d]_p \ar[r]^q & \cdot \ar[d]^g \\
        \cdot \ar[r]_f & \cdot
        }$$
    the associated composite cell
        $$\xymatrix{ \ar@{}[dr]|{\Downarrow}\cdot \ar[d]_1 \ar[r]^{p^*}|-@{|} & \cdot \ar[d]^p \ar[r]^y|-@{|} & \ar@{}[dr]|{\Downarrow} \cdot \ar[d]_q \ar[r]^{q_!}|-@{|} & \cdot \ar[d]^1 \\
        \ar@{}[dr]|{\Downarrow} \cdot \ar[d]_1 \ar[r]|-@{|} & \cdot \ar[d]^f & \ar@{}[dr]|{\Downarrow} \cdot \ar[d]_g \ar[r]|-@{|} & \cdot \ar[d]^1 \\
        \cdot \ar[r]_{f^*}|-@{|} & \cdot \ar[r]_y|-@{|} & \cdot \ar[r]_{g_!}|-@{|} & \cdot
        }$$
    is invertible.
\enddefinition

Note that precisely this composite cell appears in the proof of Lemma \ref{lemma:SufficiencyForOplaxF} where the oplax comparison cells were induced. However, Beck-Chevalley is implied by another condition involved in the definition of ``functionally complete'' (Definition \ref{define:FunctionallyComplete}) which will end up being supposed in the characterization theorem anyway.

\lemma \label{lemma:TabulatorsImpliesBeckChevalley}
    If $\D$ is an equipment with strong tabulators, then $\D_0$ has pullbacks that satisfy Beck-Chevalley. In particular, the oplax functor $F$ as in Lemma \ref{lemma:SufficiencyForOplaxF} is pseudo.
\endlemma
\proof 
    That $\D_0$ has pullbacks is proved under similar conditions in the next subsection. Proposition 5.2.3 of \cite{aleiferi} proves the entire result in detail.
\endproof

Now turn to the lax functor $G$. As in the proof of Lemma \ref{lemma:SufficiencyForLaxG} composable proarrows $p$ and $q$ induce a morphism between tabulators
    $$\xymatrix{
        \top p\otimes \top q \ar@{-->}[d]_{\gamma_{p,q}} \ar[r] & A\times C \ar@[=][d] \\
        \top (p\otimes q) \ar[r] & A\times C
    }$$
making a commutative square. If $\D_0$ is regular, then $\gamma$ is an iso if, and only if, it is a cover.

\lemma \label{lemma:TabulatorsFunctorial}
    If every relation tabulates its cokernel, then each $\gamma$ as above is an iso. In particular, the lax functor $G$ of Lemma \ref{lemma:SufficiencyForLaxG} is pseudo.
\endlemma
\proof
    The composite of tabulators in the top row of the diagram above is a relation, hence a tabulator under the assumption. Thus, the unique morphism $\gamma_{p,q}$ is an iso.
\endproof

The condition that every relation tabulates is cokernel is a powerful one. It will be discussed more in \S \ref{section:Comprehenseionscheme}. For now it is worth noting that it implies ``unit-pure.''

\lemma \label{lemma:RelationsTabulateImpliesUnitPure}
    If every relation of $\D$ tabulates its cokernel, then $\D$ is unit pure.
\endlemma
\proof
    Given morphisms $f,g\colon A\rightrightarrows B$, since the identity spans on $A$ and $B$ tabulate the corresponding identity proarrows, $f=g$ must hold by uniqueness.
\endproof

Now, the first part of the characterization result can be given.

\theorem \label{theorem:PreliminaryCharacterizationTheorem}
    Suppose that $\D_0$ is regular. The identity functor $1\colon \D_0\to \D_0$ extends to an adjoint equivalence of pseudo-functors
        \[
            F\colon \Rel(\D_0) \rightleftarrows \D \colon G
        \]
    if, and only if,
        \begin{enumerate}
            \item $y_e$ is an extension cell for each cover $e$;
            \item $\D$ is functionally complete;
            \item every relation tabulates its cokernel.
        \end{enumerate}
\endtheorem
\proof 
    The conditions are sufficient. For necessity, note that ``relations tabulate'' and each $y_e$ is an extension will imply the existence of the oplax/lax adjunction by Theorem \ref{theorem:ExistenceOplaxLaxAdjunction}. This uses Lemma \ref{lemma:RelationsTabulateImpliesUnitPure} showing that $\D$ must be unit-pure and Lemma \ref{lemma:PropertiesOfInclusionsUnitPure} showing that inclusions are then precisely the monics. So, it needs only to be seen that the lax and oplax functors are pseudo; and that these induce an adjoint equivalence.

    Strong tabulators implies Beck-Chevalley by Lemma \ref{lemma:TabulatorsImpliesBeckChevalley}. The composite cell on the left-hand side of the equation in the proof of Lemma \ref{lemma:SufficiencyForOplaxF} contains the Beck-Chevalley cell and is thus opcartesian, making the comparison cell $\phi_{R,S}$ induced there invertible, meaning that $F$ is pseudo. For $G$, that relations tabulate means that each comparison $\gamma$ is iso by Lemma \ref{lemma:TabulatorsFunctorial}, hence that each laxity comparison cells $\gamma_{p,q}$ is iso.

    The unit and counit of the induced adjunction are invertible. On the one hand, the unit $\eta_R\colon R\to \top(d^*\otimes c_!)$ is an isomorphism since $R$ is isomorphic to the tabulator of its cokernel $d^*\otimes c_!$. On the other hand, $\epsilon_p$ is iso if, and only if, $p$ is the cokernel of its tabulator. But this is precisely the condition that tabulators are ``strong" in ``functionally complete."
\endproof

\section{Comprehension Schemes} 
\label{section:Comprehenseionscheme}

The goal of this subsection is to eliminate or at least explain the unnatural condition of the previous result that every relation tabulates its cokernel. This is a sort of regularity condition on monomorphisms. It can be explained, or put in a more natural setting, by looking at generalized comprehension schemes. In fact this condition solves a problem that has not appeared yet, but will arise in the construction of a factorization system on $\D_0$ in \S \ref{section:SiteForRelations} below. This is the issue of whether tabulators provide an image factorization for $\D_0$. Specifically the issue is that of minimality, which appears to be equivalent to the statement that every relation tabulates its cokernel.

Suppose that $\D$ is a unit-pure functionally complete cartesian equipment. Taking tabulators induces a functor to subobjects
    \begin{equation} \label{equation:TabulatorAsAFunctor}
        \top\colon \D(A,B)\longrightarrow \Sub_{\D_0}(A\times B)
    \end{equation}
by the universal property of tabulators. It is well-defined by the ``monic" condition in ``functionally complete." That tabulators are strong implies that $\top$ is actually a section. But it is reasonable to ask that it is an equivalence. A couple of instructive examples are worth keeping in mind. The most immediately relevant is that of an ordinary topos $\E$. Look at this as a double category with proarrows $p\colon A\hto B$ the morphisms $p\colon A\times B\to \Omega$. Reducing to the case where $B=1$, the tabulator is simply the subobject classified by $p$ and that the function $\top$ is a bijection is precisely the statement of the universal property of the subobject classifier. Another example requires a slight stretch of the imagination. Owing to the existence of the so-called ``comprehensive factorization'' \cite{SW}, a discrete opfibration is a kind of subobject. This is precisely the approach taken in the development of 2-toposes \cite{Weber}. In the case of $\D = \prof$, the tabulator is the elements construction associated to a set-valued functor. The ``representation theorem'' for discrete opfibrations is that this has a pseudo-inverse, yielding an equivalence of categories
    \begin{equation}
        \top\colon \prof(\A,1)\longrightarrow \dopf(\A,1)
    \end{equation}
on the pattern of \ref{equation:TabulatorAsAFunctor} above. The lesson is that in each case, there can be constructed a (pseudo-)inverse to the tabulator/elements construction, yielding what could be thought of as a \textbf{comprehension scheme}. The type varies depending upon a choice of a certain monad (see \S \ref{section:ProspectusOnClassifiers} for more). In any case, think of the inverse as a \textbf{fibers construction}. In the context of topos theory, the fibers construction is the characteristic function of a subset. In the present context, this fibers construction takes the following form.

\lemma
    $\top$ as in Equation \ref{equation:TabulatorAsAFunctor} is an equivalence if, and only if, every relation $A\hto B$ is a tabulator of its cokernel.
\endlemma
\proof
    $\top$ is fully faithful by uniqueness of the arrows induced between tabulators. That relations tabulate cokernels is equivalent to the statement that $\top$ is essentially surjective. A choice of cokernels makes a strong equivalence.
\endproof

\definition \label{define:SubobjectComprehensionscheme}
    A functionally complete double category admits a \textbf{subobject comprehension scheme} if each $\top$ as in Equation \ref{equation:TabulatorAsAFunctor} is a (strong) equivalence of categories.
\enddefinition

\theorem \label{theorem:PreliminaryCharacterizationTheoremRedux}
    Suppose that $\D_0$ is regular. The identity functor $1\colon \D_0\to \D_0$ extends to an adjoint equivalence of pseudo-functors
        \[
            F\colon \Rel(\D_0) \rightleftarrows \D \colon G
        \]
    if, and only if,
        \begin{enumerate}
            \item $y_e$ is an extension cell for each cover $e$;
            \item $\D$ has a subobject comprehension scheme.
        \end{enumerate}
\endtheorem
\proof
    This is just Theorem \ref{theorem:PreliminaryCharacterizationTheorem} in light of the definition above.
\endproof

Theorem \ref{theorem:PreliminaryCharacterizationTheoremRedux} prompts consideration of conditions under which $\D_0$ is regular. This question is answered in Theorem \ref{theorem:D0isRegular} where it is shown that if $\D$ is a `double category of relations' with a subobject comprehension scheme, then $\D_0$ is regular. This leads to the main characterization result in Theorem \ref{theorem:MainCharacterizationTheorem}.

\section{Conditions for Regularity} 
\label{section:SiteForRelations}

The point of this subsection is to prove that $\D_0$ is regular under the assumption that $\D$ is a `double category of relations' with a subobject comprehension scheme. It was shown in Theorem 3.5 of \cite{CW} that $\mathbf{Maps}(\mathbb B)$, the category of maps in any `bicategory of relations', is regular and moreover that $\mathfrak{Rel}(\mathbf{Maps}(\mathbb B))\simeq\mathbb B$ as bicategories. Maps are those morphisms $p$ having a right adjoint. Here, however, the intention is to use Theorem \ref{theorem:PreliminaryCharacterizationTheoremRedux}, so $\D_0$ will be the site for relations. It is almost a category of maps in $\D$, since any ordinary morphism $f$ induces a companion and conjoint $f_!\dashv f^*$. However, an arbitrary proarrow with a right adjoint is not necessarily a companion. The point of the development is that given the other assumptions, $\D_0$ nonetheless works as a base for taking relations.

Assume that $\D$ is a `double category of relations' with a subobject comprehension scheme as in Definition \ref{define:SubobjectComprehensionscheme}. In particular, $\D$ is functionally complete and every relation tabulates it cokernel. The next two lemma construct an image factorization for any morphism of $\D_0$.

\lemma  \label{lemma:D0FactorizationExists}
    Any $f\colon A\to B$ in $\D_0$ has an image factorization $f = le$ as a cover $e$ followed by a monic $l$.
\endlemma
\proof 
    By the universal property of the tabulator there is a factorization of the cokernel of $f$
        $$\xymatrix{
            \ar@{}[dr]|{\exists \,!} A \ar[d]_e \ar[r]^y|-@{|} & A \ar[d]^e &  & \ar@{}[ddr]|{\xi}  A \ar[dd]_f \ar[r]^y|-@{|} & A \ar[dd]^f \\
            \ar@{}[dr]|{\tau} \top(f^*\otimes f_!) \ar[d]_{l} \ar[r]^y|-@{|} & \top(f^*\otimes f_!) \ar[d]^r  & = &  & \\
            B \ar[r]_{f^*\otimes f_!}|-@{|} & B &  & B \ar[r]_{f^*\otimes f_!}|-@{|} & B 
        }$$
    Since $B$ is the tabulator of $y_B$, it follows that $l=r$. Now, $\langle l, l\rangle$ is an inclusion, so by Lemmas \ref{lemma:RelationsTabulateImpliesUnitPure} and \ref{lemma:PropertiesOfInclusionsUnitPure}, it is therefore monic. Thus, $l$ itself must be monic. It needs to be seen that $e$ is a cover. For this, take any other monic $m\colon M\to B$ through which $f$ factors by say $p\colon A\to M$. The cokernel of $f$, since it is an opcartesian cell, factors through the cokernel of $m$. Therefore, since $M$ tabulates the cokernel of $\langle m,m\rangle$, $l$ factors through $m$ uniquely. In other words, $l$ is the smallest subobject through which $f$ factors, meaning that $f=le$ is an image factorization.
\endproof

\remark
    It is worth pausing here to explain the importance of the assumption that relations tabulate. This was used above to prove essentially that $e$ is extremal. In each of \cite{CW}, \cite{FS}, and \cite{Elephant}, there is only one level of arrow, so that the corresponding map is extremal is a consequence of the fact that it turns out to be both simple and entire. However, this same move cannot be made in the present context of double categories. For the induced arrow to the tabulator does not provably have an inverse from properties of its companion and conjoint. Precisely what is needed is the ``fibers construction'' coming with the subobject comprehension scheme.
\endremark

Now that $\D_0$ has its factorization system, for regularity, it must be shown that covers are pullback-stable. First it needs to be seen that pullbacks exist. The proof follows closely that of Proposition 3.2.7 in \cite{Elephant}. See also Proposition 5.2.3 of \cite{aleiferi}.

\lemma \label{lemma:pullbacksExist}
    $\D_0$ has all pullbacks.
\endlemma
\proof 
    Take a corner diagram $h\colon A\to C \leftarrow B\colon e$. Take the tabulator of the restriction as in the diagram
        $$\xymatrix{
            \ar@{}[dr]|{\tau} \top (h_!\otimes e^*) \ar[d]_d \ar[r]^y|-@{|} & \top (h_!\otimes e^*) \ar[d]^c \\
            \ar@{}[dr]|{\xi} A \ar[d]_h \ar[r]^{h_!\otimes e^*}|-@{|} & B \ar[d]^e \\
            C \ar[r]_y|-@{|} & C
        }$$
    The arrows $d$ and $c$ now complete the pullback square
        $$\xymatrix{ 
            P \ar[d]_d \ar[r]^c & B \ar[d]^e \\
            A \ar[r]_h & C
        }$$
    The square commutes because $C$ is the tabulator of $y_C$. The universal property for the pullback follows by the universal property of the tabulator.
\endproof

Now it can be seen that $\D_0$ is regular. Pullback-stability of covers is all that remains. Recall that ``functionally complete" implies the Beck-Chevalley condition by Lemma \ref{lemma:TabulatorsImpliesBeckChevalley}.

\theorem \label{theorem:D0isRegular}
    If $\D$ is cartesian, then $\D_0$ is a regular category, hence $\Rel(\D_0)$ is well-defined.
\endtheorem
\proof 
    Since $\D$ is cartesian, $\D_0$ has finite products. Existence of pullbacks was proved above in Lemma \ref{lemma:pullbacksExist}. It needs only to be seen that covers are stable under pullback. For this let $e\colon B\twoheadrightarrow C$ denote a cover. The pullback square
        $$\xymatrix{ 
            P \ar[d]_d \ar[r]^c & B \ar@{->>}[d]^e \\
            A \ar[r]_h & C
        }$$
    is formed using the tabulator of the restriction $h_!\otimes e^*$ as in Lemma \ref{lemma:pullbacksExist}. 
    On the one hand, since $h_!\otimes h^*$ is a restriction, there is a canonical cell $y_A\Rightarrow h_!\otimes h^*$. To see that $d$ is a cover, calculate that
        \begin{align} 
            y_A & = y\wedge (h_!\otimes h^*)\notag \\
              & = y\wedge (h_!\otimes e^*\otimes e_!\otimes h^*) \qquad &\text{($e$ is a cover)} \notag\\
              & = y\wedge (d^*\otimes c_!\otimes c^*\otimes d_!) \qquad &\text{(Beck-Chevalley \ref{define:BeckChevalley})} \notag\\
              & \leq d^*\otimes(d_!\wedge c_!\otimes c^*\otimes d_!) \qquad &\text{(Modular Law \ref{equation:ModularLaw})} \notag\\
              & \leq d^*\otimes d_!. \qquad \notag
        \end{align}
    Since there is always a canonical map $d^*\otimes d_!\Rightarrow y$, this shows by the hypothesis that $y$ is fully faithful that $d^*\otimes d_!\cong y$ holds, meaning that $d$ is indeed a cover, proving that $\D_0$ is regular.
\endproof

\remark
    A shorter proof of Theorem \ref{theorem:D0isRegular} might simply quote the criteria of Theorem 4.4.4 in \cite{BJ}. This result says that a finitely-complete category $\E$ is regular if, and only if, the subobject fibration $\mathrm{Sub}(\E)\to \E$ has coproducts satisfying a corresponding Frobenius identity. This seems a rather easy criterion to use since under the hypotheses of the subsection, inclusions are monic and the Modular Law \ref{equation:ModularLaw} is very close to what is meant by ``Frobenius'' in the reference. However, without proving the equivalence of the conditions in the quoted result, it is not clear that this approach would be of considerable benefit. The space-saving otherwise would result in a loss of explicitness in the constructions whereas on the other hand proving the theorem would involve for the most part reproducing the work done here already. 
\endremark

\section{Characterization Theorem}
\label{section:CharacterizationTheorem}

Now, the main result of the paper can be given. The goal is to use Theorem \ref{theorem:PreliminaryCharacterizationTheoremRedux}. For this note the following result, saying, essentially, that ``covers are covers.''

\lemma \label{lemma:coversarecovers}
    Let $\D$ denote a `double category of relations' with a subobject comprehension scheme. Any cover in the regular category $\D_0$ is then one in $\D$.
\endlemma
\proof
    First show that $e$ as in the proof of Lemma \ref{lemma:D0FactorizationExists} is a cover in the equipment structure. Note that $e_!= f_!\otimes l^*$ holds and dually $e^*= l_!\otimes f^*$ holds too, since in each case either side computes the same restriction. Now, $f$ and $l$ have the same cokernel because tabulators are strong. This means that 
        \[ 
            e^*\otimes e_!= l_!\otimes f^*\otimes f_!\otimes l^*= l_!\otimes l^* \otimes l_!\otimes l^* = y\otimes y= y
        \]
    holds since $l$ is monic, hence an inclusion, proving that $e$ is a cover in $\D$. This means that every morphism factors uniquely as a cover in the equipment structure followed by a monic arrow in $\D_0$. In particular, every cover in $\D_0$ is one in $\D$. For in this case $l$ is invertible, hence a cover too, meaning that $y\cong l^*\otimes l_!\cong f^*\otimes f_!$ must hold since $f$ and $l$ have the same cokernel.
\endproof

Now, the main result of the paper, characterizing double categories of relations as `double categories of relations' with a subobject comprehension scheme.

\theorem \label{theorem:MainCharacterizationTheorem}
    If $\D$ is a `double category of relations' with a subobject comprehension scheme then the identity functor $1\colon \D_0\to \D_0$ extends to an adjoint equivalence 
        \[ 
            \Rel(\D_0)\simeq \D.
        \]
    In short, any `double category of relations' with a subobject comprehension scheme is equivalent to a double category $\Rel(\E)$ for some regular category $\E$.
\endtheorem
\proof 
    Theorem \ref{theorem:D0isRegular} shows that $\D_0$ is regular. That the identity functor on $\D_0$ then extends to an equivalence is then Theorem \ref{theorem:PreliminaryCharacterizationTheoremRedux} by way of Lemma \ref{lemma:coversarecovers}.
\endproof

Now that the characterization has been given, there are several immediate results. Recall that allegories are defined in such a way as to possess an \textbf{anti-involution operator}. This is integral to the definition of the Modular Law. In an allegory of relations, the involution is just to take the opposite relation. Double categories of relations have a derived operation coming from the existence of tabulators and other exactness conditions. For the second statement of the following cf. A3.2.3 of \cite{Elephant}.

\proposition \label{proposition:InvolutionOperatorExists}
    A double category of relations has an anti-involution operator, that is, an operation on proarrows $R\mapsto R^\circ$ such that $y^\circ = y$ and $(R^\circ)^\circ= R$ both hold. For a proarrow of the form $p_!$, the involution is the same as the right adjoint $(p_!)^\circ= p^*$.
\endproposition
\proof
    Given a proarrow $R\colon A\hto B$, take $R^\circ$ to be the cokernel of the opposite of the tabulator of $R$. For the second statement, if the tabulator of $p_!$ has legs $l$ and $r$, then by lifting and extension properties, the cokernel $r^*\otimes l_!$ is right adjoint to $p_!$ in the proarrow bicategory of $\D$ and so must be isomorphic to $p^*$ canonically.
\endproof

Double categories of relations are also regular as equipments. This is an immediate corollary since every double category $\Rel(\E)$ is regular \cite{Schultz}. It is, however, worth spelling out in somewhat more detail what this means.

Recall first a few preliminary definitions. For those of a monoid, bimodule and their homomorphisms one can see the reference, \S 2 of \cite{FrameworkGenMulticats}, or \cite{discdblefibs}. A monoid in $\D$ is \textbf{effective} if it is the kernel of some ordinary morphism. An \textbf{embedding} of a monoid $m\colon A\hto A$ into an object $X$ is a monoid homomorphism $m\to X$ from $m$ to the trivial monoid on $X$. The \textbf{collapse} of a monoid is a universal embedding. Likewise a \textbf{bimodule collapse} is a universal bimodule embedding. The collapse of a monoid is \textbf{normal} if it presents the bimodule collapse of the trivial bimodule on $A$. A morphism $f$ is a \textbf{regular cover} if its kernel is a normal collapse cell.

\definition[Definition 4.7, \cite{Schultz}]
    A double category $\D$ is \textbf{regular} if
        \begin{enumerate}
            \item every effective monoid has a normal collapse;
            \item every restriction cell
                    $$
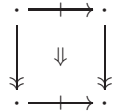
\xymatrix{
                        \ar@{}[dr]|{\Downarrow} \cdot \ar@{->>}[d] \ar[r]|-@{|} & \cdot \ar@{->>}[d] \\
                        \cdot \ar[r]|-@{|} & \cdot
                    }$$
                is a bimodule collapse cell.
        \end{enumerate}
\enddefinition

The first condition is the analogue of the condition that every kernel has a coequalizer; the second is the condition that regular epimorphisms are pullback-stable.

\corollary \label{corollary:DoubleCategoriesOfRelationsAreRegular}
    Any double category of relations with a subobject comprehension scheme is regular as an equipment.
\endcorollary
\proof
    Any such double category is up to equivalence of the form $\Rel(\E)$ for some regular category $\E$ by Theorem \ref{theorem:MainCharacterizationTheorem}, which is regular by Proposition 4.8 in the reference. 
\endproof

In a regular equipment, the canonical factorization of a morphism $f\colon A\to B$ developed in \S 4.2 of \cite{Schultz} is given by taking the collapse of the kernel of $f$. This is meant to mimic the factorization system in ordinary regular categories given by taking the coequalizer of the kernel of a given morphism. The approach in the development of double categories of relations here has been a dual construction, namely, taking the tabulator of the cokernel of a morphism. So, on the one hand, for regular equipments, the factorization is a quotient of a kernel, for double categories of relations the factorization is a subobject of a cokernel. However, owing to the further exactness conditions in a double category of relations, these two factorizations coincide.

\proposition
    The factorization of a morphism in a double category of relations viewed as a regular double category coincides with the factorization produced in Lemma \ref{lemma:D0FactorizationExists} above.
\endproposition
\proof
    The factorization as a regular double category is as a regular cover followed by an inclusion. However, inclusions here have the same definition as in the reference. The condition that there is a fibers construction as in Definition \ref{define:BeckChevalley} means that any inclusion is a tabulator, hence that the two factorizations coincide.
\endproof

\section{Division and Powers}
\label{section:DivisionPowers}

Recall \cite[\S I,7]{FS} that a \textbf{logos} is a regular category $\E$ such that each $\Sub(A)$ is a lattice and every pullback functor $f^*\colon \Sub(B)\to \Sub(A)$ has a right adjoint $\forall_f$. A logos is called a ``Heyting category'' in \cite{Elephant}. Thus, for a logos, each pullback functor $f^*$ has both a left and a right adjoint with the left adjoint $\exists_f$ given by taking images. The following is the double-categorical version of a division allegory, that is, an allegory equipped locally with certain division operators (cf. \cite[\S II,2.31]{FS} or \cite[\S A3.4]{Elephant}).

\definition
    A double category of relations $\D$ has \textbf{division} if for any proarrow $p\colon A\to B$ and any object $C$, the functor 
        \[
          p\otimes (-)\colon \D(B,C) \to \D_1(A,C)
        \]
    given by precomposition with $p$ has a right adjoint. Denote the right adjoint by $(-)/p$.
\enddefinition

Double categories of relations with division are related to logoi in the following way. This gives an appropriate version of \S II,2.32 of \cite{FS} which says that an allegory is a division allegory if, and only if, its underlying category of maps is a logos.

\theorem
    A double category of relations $\D$ has division if, and only if, each pullback functor $f^*$ between subobject posets has a right adjoint. Consequently, a double category of relations with a subobject comprehension scheme and division is of the form $\Rel(\E)$ with $\E$ a logos if, and only if, each subobject poset is a lattice.
\endtheorem
\proof
    Theorem \ref{theorem:MainCharacterizationTheorem} shows that any double category of relations is of the form $\Rel(\E)$ for some regular category $\E$. So, on the one hand, if $\Rel(\E)$ has division, then identifying $\Rel(C)_1(X,1) \simeq \Sub(X)$ for any object $X$, each pullback functor between posets has a right adjoint. On the other hand, suppose that each pullback functor between subobject posets has a right adjoint. Fix a relation $R\colon A\hto B$ with legs $l$ and $r$. Both functors on the bottom row of the commutative diagram 
        $$\xymatrix{
            \Rel(C)_1(B,C) \ar[d]_\simeq \ar[r] & \Rel(C)_1(R,C) \ar[d]_{\simeq} \ar[r] & \Rel(C)_1(A,C) \ar[d]^\simeq \\
            \Sub(B\times C) \ar[r]_{(r\times 1)^*}  & \Sub(R\times C) \ar[r]_{\exists_{l\times 1}} & \Sub(A\times C) 
        }$$
    have right adjoints. The top row is first restriction along $r$ and then extension along $l$. Thus, up to the equivalences above, the right adjoint is $\forall_{r\times 1}\circ (l\times 1)^*$. The last statement follows by the foregoing result and the definition of a logos recalled above.
\endproof

\remark
    Notice that in the background the identification $\D_1(A,B)\simeq \Sub(A\times B)$ given by the tabulator in Definition \ref{define:SubobjectComprehensionscheme} is actually doing most of the work in the proof of the theorem. 
\endremark

Now recall that a \textbf{topos} is a finitely complete category $\E$ with \textbf{power objects}, namely, special objects $PA$ for each object $A$ and special monomorphisms $\in_A\rightarrowtail PA\times A$ that are suitably universal among subobjects of the form $S\rightarrowtail X\times A$. Universality is expressed by a pullback condition. This is equivalent to the usual standard definition \cite[\S A2.1]{Elephant} as a cartesian closed category with a subobject classifier. An allegory is a \textbf{power allegory} if each object $A$ has a map $PA\to A$ satisfying a couple of technical conditions, namely, equations stating that the morphism $PA\to A$ satisfies extensionality and comprehension. The result of \S II,2.414 of \cite{FS} is that any tabular power allegory is of the form $\mathbf{Rel}(\E)$ where $\E$ is a topos. See also Corollary A3.4.7 of \cite{Elephant}. In the present context, the universal property of power objects is especially easy to state using just the equipment axioms. Recall that restrictions in an equipment should be thought of as pullbacks.

\definition
    An equipment $\D$ has \textbf{powers} if each object $A$ is equipped with a special proarrow $\in_A \colon PA\hto A$ such that for any $R\colon X\hto A$ there is a unique morphism $X\to PA$ yielding a restriction cell
        $$\xymatrix{
            \ar@{}[dr]|{\Downarrow} X \ar[d] \ar[r]|-@{|} & A \ar@{=}[d] \\
            PA \ar[r]_{\in_A}|-@{|} & A
        }$$
    Equivalently, $\D$ has powers if $\D(-,A)\colon \D_0^{op}\to \set$ is representable.
\enddefinition

\theorem
    Any double category of relations with a subobject comprehension scheme and powers is of the form $\Rel(\E)$ for $\E$ a topos.
\endtheorem
\proof
    Since restrictions in $\Rel(\E)$ are computed by pullback, the universal property of $\in_A$ above is precisely that of the universal subobject $\in_A\rightarrowtail PA\times A$ in the equivalent definition of a topos.
\endproof

\section{Prospectus}
\label{section:Prospectus}

Let us end with some comments on future and ongoing work. In particular, there are many potential applications, mostly within category theory, but possibly in other areas. Several were mentioned in the introduction. Here a few can be discussed in somewhat more detail. 

\subsection{Double-Categorical Semantics}

It was remarked in \S \ref{section:DoubleCategoriesOfRelations} that the horizontal bicategory $\mathcal H(\D)$ of any cartesian equipment $\D$ might be regarded as a generalized compact closed category. This is essentially a consequence of Proposition \ref{prop:HorizontalBicatIsCartesianBicat} showing that $\mathcal H(\D)$ is a cartesian bicategory. More precisely, a compact closed category \cite{KellyLaplaza} is a symmetric monoidal category with a certain dualization operation. As such, compact closed categories are ``$*$-autonomous'' in the sense of \cite{BarrAutonomous}. It is well-known that $*$-autonomous categories provide a categorical semantics for linear logic \cite{Seely}, and it seems \emph{prima facie} possible to import much of the development to the context of cartesian equipments. Since resources are interpreted as objects in a $*$-autonomous category, in the double categorical context, they would be interpreted as proarrows, with cells providing interpretations of deductions. This allows the possibility of introducing context-dependency for deployment of resources and of an underlying type theory that is modeled by the objects and ordinary arrows of the cartesian double category. The ``logic'' part of linear logic would then be interpreted by whatever extra structure and local connectives were asked for in $\D_1$.

This is part of a broader program, that owing to the fact that any equipment $\D$, being at least a fibration, has both an internal type theory given by $\D_0$ and an internal logic of types and terms given by $\D_1$. This makes sense as proarrows should be interpreted as relations or their generalizations such as spans or honest profunctors. In ordinary first-order logic, predicate symbols and thus formulas are interpreted categorically by subobjects, which technically speaking live in a separate poset. The virtue of modeling type theories and their logics in suitably structured double categories is that the proarrows and various local connectives (i.e. local products) have enough structure to deal with the type-theoretic and logical aspects in the same structure. For example, the interpretation of regular logic in bicategories of relations \cite{patterson}, \cite{fongspivak}, which in the former is limited to predicate symbols only, can be done for full regular logic with predicate and function symbols by working in a cartesian double category instead. 

\subsection{Monoidal Fibrations}

In \S 14 of \cite{FramedBicats}, it is shown that under certain conditions, every monoidal bifibration gives rise to an equipment with some extra structure. However, it appears that the definition of a cartesian equipment, as presented in \cite{aleiferi}, is needed for a sort of inverse construction taking a cartesian equipment to a monoidal bifibration. Under such a hypothetical correspondence, it is of interest to see which monoidal bifibrations correspond to double categories of relations. The conjecture is that these will be closely related to regular fibrations and subobject fibrations as in \S 4.2 and \S 4.4 of \cite{BJ}. If this is the case, it is another point in favor of the view of the close connection between type theories and suitably structured double categories.

\subsection{Classification Schemes for Double Toposes}
\label{section:ProspectusOnClassifiers}

As justification for the name ``subobject comprehension scheme'' adopted in Definition \ref{define:SubobjectComprehensionscheme}, recall some of the development of hyperdoctrines from \cite{hyper}. An \textbf{elementary existential doctrine} is a pseudo-functor on a cartesian category $P\colon \E^{op}\to \cat$ such that each substitution functor $f^*\colon PB\to PA$ has a left adjoint $\Sigma_f$. Suppose that each category $PA$ has a terminal object. There is then a natural functor $\E/B\to PB$ taking a morphism $f\colon A\to B$ to $\Sigma_f(1)$. If each such functor has a right adjoint $\lbrace -\rbrace$, then $P$ is a \textbf{comprehension scheme}. In the case that $\D$ is a cartesian equipment with tabulators, the hyperdoctrine
    \[
      \D(-,1)\colon \D_0^{op}\to\cat  
    \]
is an elementary existential doctrine with extension providing the left adjoint to restriction as substitution. The comprehension scheme is then given by tabulators
    \[
      \Sigma_{(-)}1\colon \D_0/B \leftrightarrows \D(B,1)\colon \top.
    \]
If $\D_0$ is regular and tabulators are monic, tabulators factor through the subobject poset as in 
    $$\xymatrix{
        \D_0/B \ar@<-.5ex>[rr]_{\Sigma} \ar@<-.5ex>[dr]_\sigma & & \ar@<-.5ex>[ll]_\top \D(B,1) \ar[dl]^\top \\
        & \ar@<-.5ex>[ul]_i \Sub(B)
    }$$
where $\mathrm{\sigma}$ is left adjoint to the inclusion of subobjects. What the ``subobject comprehension schemes'' of \S \ref{section:Comprehenseionscheme} axiomatizes is that the left adjoint $\Sigma$, when restricted to subobjects, results in an equivalence making the other triangle above commute and identifying $\D(B,1)$ as a reflective subcategory of the slice. 

Although strictly speaking \cite{hyper} \textit{defines} the comprehension scheme as being the mere presence of the right adjoint in the top row of the diagram above, the present view is that the structure given by this right adjoint should be taken into account. In the example under discussion, the tabulator gives a monic arrow which represents a subobject. In the other example, namely, that of $\D = \prof$, tabulators are not monic, but instead are discrete opfibrations. There results a situation somewhat like that above
    $$\xymatrix{
        \cat/\B \ar@<-.5ex>[rr]_{\Sigma} \ar@<-.5ex>[dr]_\sigma & & \ar@<-.5ex>[ll]_\elt \prof(\B,1) \ar[dl]^\elt \\
        & \ar@<-.5ex>[ul]_i \dopf(\B)
    }$$
where now $\sigma$ gives the free discrete opfibration on a given functor over $\B$. Again there is a fibers construction assigning to every discrete opfibration a set-valued functor, that is, a profunctor $\B\to 1$. This results in the well-established equivalence between discrete opfibrations and set-valued functors. This fibers construction gives the analogue of the characteristic function from ordinary topos theory. The reason for treating this part of the triangle as the distinctive aspect of the comprehension scheme is simply that fibers do not vary functorially unless the construction starts with a discrete opfibration. Hence the language ``subobject comprehension scheme'' or ``discrete opfibration comprehension scheme'' clearly identifying the properties of the projection morphism coming with the right adjoint.

Now, discrete opfibrations are not reflective in $\cat/\B$, but there is a commonality with the previous example in that both categories $\Sub(B)$ and $\dopf(\B)$ are monadic over the respective slice categories. In particular, discrete opfibrations are algebras for a ``pull-push monad" on $\cat/\B$ described for example in \S 2.2-2.3 of \cite{TT}. What is envisioned is that $\D$ is a cartesian equipment with tabulators, each slice $\D_0/B$ is equipped with a monad $T$ and tabulators factor through $T$-algebras as in
    $$\xymatrix{
        \D_0/B \ar@<-.5ex>[rr]_{\Sigma} \ar@<-.5ex>[dr]_\sigma & & \ar@<-.5ex>[ll]_\top \D(B,1). \ar[dl]^\top \\
        & \ar@<-.5ex>[ul]_i T\alg
    }$$
For this to be a $T$-\textbf{comprehension scheme} required would be a fibers construction from $T$-algebras back to proarrows resulting in an equivalence of categories. Such a template recovers the examples discussed so far. Toposes would be a special case where there is a representing object, namely, the subobject classifier $\Omega$ for the original hyperdoctrine. Such structures could figure prominently in double-categorical interpretations of higher-order type theory (cf. the ``comprehension categories'' in Ch. 10 of \cite{BJ}).

This approach could avoid the issue of ``admissibility'' in the development of 2-toposes \cite{Weber} and possibly Yoneda structures generally \cite{SW} where size issues prevent the development of a genuine pseudo-inverse to the elements construction (cf. Definition 4.1 of \cite{Weber} where the functor that might be expected to be an equivalence is merely fully faithful). This is because of the requirement that $\set$ lives in an enlarged $\cat$ and so there is no guarantee that an arbitrary discrete opfibration has small fibers. The approach being suggested here is that $\set$ is \textit{supposed to be} the representing object for the hyperdoctrine, but $\set$ is ``hidden behind the proarrows" coming with the double category structure. In this way, the fibers construction makes sense yielding the genuine equivalence without requiring the existence of a representing object for the hyperdoctrine. So, in this sense, it is plausible that 2-toposes are rather 2-categorical fragments of ``double toposes,'' which should be at least cartesian equipments with $T$-comprehension schemes for certain well-chosen monads $T$. This is further supported by the fact that a nice enough 2-topos ends up supporting a Yoneda structure \cite{Weber} which in any case behaves very much like the proarrows of an equipment \cite{FramedBicats}.

\refs

\bibitem[Aleiferi, 2018]{aleiferi} Evangelia Aleiferi. ``Cartesian Double Categories with an Emphasis on Characterizing Spans," PhD Thesis, Dalhousie University, 2018, arXiv:1809.06940.

\bibitem[Barr, 1979]{BarrAutonomous} Michael Barr. *-Autonomous Categories. Volume 752 of \textit{Lecture Notes in Mathematics}. Springer, Berlin, 1979. 

\bibitem[Bonchi et. al., 2017]{BFS} Filippo Bonchi, Dusko Pavlovic, and Pawel Sobocinski. ``Functorial Semantics for Relational Theories.'' arXiv:1711.08699, 2017.

\bibitem[Carboni \& Walters, 1987]{CW} A. Carboni and R.F.C. Walters. ``Cartesian Bicategories I." \textit{Journal of Pure and Applied Algebra}, 49(1), pp. 11-32, 1987.

\bibitem[Codd, 1972]{Codd} E.F. Codd. ``Relational Completeness of Data Base Sublanguages." In: \textit{Database Systems}. Prentice-Hall, pp. 65-98, 1972.

\bibitem[Cruttwell \& Shulman, 2010]{FrameworkGenMulticats} G.S.H. Cruttwell and M.A. Shulman. ``A Unified Framework for Generalized Multicategories." \textit{Theory and Applications of Categories}, 24(21), pp. 580-655, 2010.

\bibitem[Di Liberti et. at., 2021]{LLNSob} Ivan Di Liberti, Fosco Loregian, Chad Nester, and Pawe{\l} Soboci\'nski. ``Functorial Semantics for Partial Theories.'' \textit{Proceedings of the ACM on Programming Languages}, 5(57), pp. 1-28, 2021.

\bibitem[Ehresmann, 1963]{Ehresmann} C. Ehresmann. ``Cat\'egories et structures." Dunod, Paris, 1963.

\bibitem[Fong \& Spivak, 2019]{fongspivak} Brendan Fong \& David Spivak. ``Regular and Relational Categories: Revisiting `Cartesian Bicategories I.''' arXiv:1909.00069.

\bibitem[Freyd \& Scedrov]{FS} Peter J. Freyd and Andre Scedrov. \textit{Categories, Allegories}. North Holland, 1990.

\bibitem [Grandis \& Par\'e, 1999]{GP} M. Grandis and R. Par\'e. ``Limits in Double Categories." \textit{Cahiers de Topologie et G\'eom\'etrie Diff\'erentielle Cat\'egoriques}, 40(3), pp. 162-220, 1999.

\bibitem [Grandis \& Par\'e, 2004]{GP2} M. Grandis and R. Par\'e. ``Adjoints for Double Categories." \textit{Cahiers de Topologie et G\'eom\'etrie Diff\'erentielle Cat\'egoriques}, 45(3), pp. 193-240, 2004.

\bibitem[Jacobs, 1999]{BJ} Bart Jacobs. \textit{Categorical Logic and Type Theory}. Volume 141 of \textit{Studies in Logic and the Foundations of Mathematics}, Elsevier, 1999.

\bibitem[Johnstone, 1977]{TT} P.T. Johnstone. \textit{Topos Theory}. Academic Press, New York, 1977.

\bibitem[Johnstone, 2001]{Elephant} P.T. Johnstone. \textit{Sketches of an Elephant: A Topos Theory Compendium, Volume 1}. Volume 43 of \textit{Oxford Logic Guides}. Clarendon Press, London, 2001.

\bibitem[Kelly \& Laplaza]{KellyLaplaza} G.M. Kelly and M.L. Laplaza. ``Coherence for Compact Closed Categories.'' \textit{Journal of Pure and Applied Algebra}, 19, pp. 193-213, 1980.

\bibitem[Kent \& Spivak, 2012]{SpivakKent} Robert E. Kent and David I. Spivak. ``Ologs: A Categorical Framework for Knowledge Representation.'' \textit{PLoS ONE}, 7(1), 2012.

\bibitem[Knijnenburg \& Nordemann, 1994]{PetrusAndFrank} Petrus Knijnenburg and Frank Nordemann. ``Two Categories of Relations.'' \textit{Technical Report, Leiden University Department of Computer Science}, 94(32), 1994.

\bibitem[Koudenburg, 2015]{roalddoubledim} Seerp Roald Koudenburg. ``A Double-Dimensional Approach to Formal Category Theory.'' arXiv:1511.04070, 2015.

\bibitem[Koudenburg, 2019]{roaldaugmented} Seerp Roald Koudenburg. ``Augmented Virtual Double Categories.'' \textit{Theory and Applications of Categories}, 35(10), pp. 261-325, 2019.

\bibitem[Lambert, 2021]{discdblefibs} Michael Lambert. ``Discrete Double Fibrations.'' \textit{Theory and Applications of Categories}, 37(22), 2021, pp. 671-708.

\bibitem[Lambert, 2022]{lamberttoposblog} Michael Lambert. ``Data Operations are Functorial Semantics.'' \textit{The Topos Lab}, 2022, https://topos.site/blog/2022/09/data-operations-are-functorial-semantics/.

\bibitem[Lawvere, 1970]{hyper} William F. Lawvere. ``Equality in Hyperdoctrines and Comprehension Schema as an Adjoint Functor.'' Volume XVII of \textit{Proceedings of the AMS Symposium on Pure Mathematics}, pp. 1-14, 1970.

\bibitem [MacLane, 1998]{MacLane} S. MacLane. \textit{Category Theory for the Working Mathematician}. Volume 5 of \textit{Graduate Texts in Mathematics}, Springer, Berlin, 1998.

\bibitem[Niefield, 2012]{Niefield} Susan Niefield. ``Span, Cospan, and Other Double Categories.'' \textit{Theory and Applications of Categories} 26(26), pp.729-742, 2012.

\bibitem[Patterson, 2017]{patterson} Evan Patterson. ``Knowledge Representation in Bicategories of Relations." arXiv:1706.00526, 2017.

\bibitem[Rosebrugh \& Wood, 1992]{RW} Robert Rosebrugh and Richard Wood. ``Relational Databases and Indexed Categories.'' \textit{Category Theory 1991: Proceedings of the 1991 Summer Category Theory Meeting, Montreal, Canada}. Volume 13 of \textit{Canadian Mathematical Society Conference Proceedings}, pp. 390-407, 1992.

\bibitem[Schultz, 2015]{Schultz} Patrick Schultz. ``Regular and Exact (Virtual) Double Categories.'' arXiv:1505.00712, 2015.

\bibitem [Shulman, 2008]{FramedBicats} M. Shulman. ``Framed bicategories and monoidal fibrations." \textit{Theory and Applications of Categories}, 20(18), pp. 650-738, 2008.

\bibitem[Selinger, 1999]{SelingerAsynchronous} Peter Selinger. ``Categorical Structure of Asynchrony.'' \textit{Proceedings of the 15th Conference on the Mathematical Foundations of Programming Semantics, MFPS 1999, New Orleans}. Volume 20 of \textit{Electronic Notes in Theoretical Computer Science}, pp. 158-181, 1999.

\bibitem[Seely, 1989]{Seely} R.A.G. Selly. ``Linear Logic, {*-}Autonomous Categories and Cofree Coalgebras.'' In. John W. Gray ed. \textit{Category Theory in Computer Science and Logic}, Volume 92 of \textit{Contemporary Mathematics}, American Mathematical Society, pp. 371-383, 1989.

\bibitem[Street \& Walters, 1973]{SW} Ross Street and R.C.F. Walters. ``The Comprehensive Factorization of a Functor.'' \textit{Bulletin of the American Mathematical Society}, 79(5), pp. 936-941, 1973.

\bibitem[Weber, 2007]{Weber} Mark Weber. ``Yoneda Structures from 2-Toposes.'' \textit{Applied Categorical Structures}, 15, pp. 259-323, 2007.

\bibitem [Wood, 1982]{WoodProI} R. Wood. ``Abstract proarrows I." \textit{Cahiers de Topologie et G\'eom\'etrie Diff\'erentielle Cat\'egoriques}, 23(3), pp. 279-290, 1982.

\bibitem [Wood, 1985]{WoodProII} R. Wood. ``Proarrows II." \textit{Cahiers de Topologie et G\'eom\'etrie Diff\'erentielle Cat\'egoriques}, 26(2), pp. 135-168, 1985.

\endrefs

\end{document}